\documentclass[aos,MSNbibl,nameyear,seceqn,dvips]{arximspdf}
\usepackage{multirow}
\usepackage{graphicx}

\doi{10.1214/12-AOS1016} 
\volume{40}
\issue{3}
\pubyear{2012}
\firstpage{1609}
\lastpage{1636}

\makeatletter
\newcommand{\rrVert}{\Vert}
\newcommand{\rrvert}{\vert}
\newcommand{\llVert}{\Vert}
\newcommand{\llvert}{\vert}
\newtheorem{thmm}{Theorem}[section]
\newproclaim{rem}{Remark}[section]
\makeatother

\begin{document}
\begin{frontmatter}

\title{Partially monotone tensor spline estimation of the joint
distribution function with bivariate current~status data}
\runtitle{Sieve estimation with bivariate current status data}

\begin{aug}
\author[A]{\fnms{Yuan} \snm{Wu}\ead[label=e1]{y6wu@ucsd.edu}}
\and
\author[B]{\fnms{Ying} \snm{Zhang}\corref{}\ead[label=e2]{ying-j-zhang@uiowa.edu}}
\runauthor{Y. Wu and Y. Zhang}
\affiliation{University of California San Diego and University of Iowa}
\address[A]{Division of Biomedical Informatics\\
University of California, San Diego\\
La Jolla, California 92093\\
USA\\
\printead{e1}}

\address[B]{Department of Biostatistics\\
University of Iowa\\
Iowa City, Iowa 52242\\
USA\\
\printead{e2}}

\end{aug}

\received{\smonth{12} \syear{2011}}
\revised{\smonth{5} \syear{2012}}

%
\begin{abstract}
The analysis of the joint cumulative distribution function (CDF)
with bivariate event time data is a challenging problem both
theoretically and numerically. This paper develops a tensor
spline-based sieve maximum likelihood estimation method to estimate
the joint CDF with bivariate current status data. The $I$-splines
are used to approximate the joint CDF in order to simplify the
numerical computation of a constrained maximum likelihood estimation
problem. The generalized gradient projection algorithm is used to
compute the constrained optimization problem. Based on the
properties of $B$-spline basis functions it is shown that the
proposed tensor spline-based nonparametric sieve maximum likelihood
estimator is {consistent} with a~rate of convergence {potentially}
better than $n^{1/3}$ under some mild regularity conditions. The
simulation studies with moderate sample sizes are carried out to
demonstrate that the finite sample performance of the proposed
estimator is generally satisfactory.
\end{abstract}

%


\end{frontmatter}

\section{Introduction}\label{sec1}

In some applications, observation of random event\break time~$T$ is
restricted to the knowledge of whether or not $T$ exceeds a random
monitoring time $C$. This type of data is known as current status
data and sometimes referred to as interval censoring case 1
[\citet{GroWel92}]. Current status data arise naturally
in many applications; see, for example, the animal tumorigenicity
experiments by \citet{HoeHEW72} and \citet{FinWol85}; the social demographic studies of the distribution of the age
at weaning by \citet{DiaMcDSha86}, \citet{DiaMcD91} and \citet{Gru93}; and the studies of human
immunodeficiency virus (HIV) and acquired immunodeficiency syndrome
(AIDS) by \citet{ShiJew92} and \citet{JewMalVit94}.

The univariate current status data {have} been thoroughly studied in
the statistical literature. \citet{GroWel92} and \citet{HuaWel95} studied the asymptotic properties of the
nonparametrc maximum likelihood estimator (NPMLE) of the CDF with
current status data. \citet{Hua96} considered Cox proportional
hazards model with current status data and showed that the maximum
likelihood estimator (MLE) of the regression parameter is
asymptotically normal with $\sqrt{n}$ convergence rate, even through
the MLE of the baseline cumulative hazard function only converges at
$n^{1/3}$ rate.


Bivariate event time data occur in many applications as well. For
example, in an Australian twin study [\citet{DufMarMat90}], the researchers were interested in times to a certain event
such as a disease or a disease-related symptom in both twins. NPMLE
of the joint CDF of the correlated event times with bivariate right
censored data was studied by \citet{Dab88}, \citet{PreCai92}, \citet{Pru}, \citet{van96}
and \citet{QuavanRob06}. As an alternative, \citet{Koo98} developed a
tensor spline estimation of the logarithm of joint density function
with bivariate right censored data. However, asymptotic properties
of Kooperberg's estimate are unknown. \citet{ShiLou95} proposed
a two-stage semiparametric estimation procedure to study the joint
CDF with bivariate right censored data, in which the joint
distribution of the two event times is assumed to follow a bivariate
Copula model [\citet{Nel06}].


For bivariate interval censored data, the {conventional} NPMLE {was
originally} studied by \citet{BetFin99} {and
followed by} \citet{WonYu99}, \citet{GenVan01}, \citet{Son01} and \citet{Maa05}.
{A typical numerical algorithm for
computing the NPMLE constitutes two steps [\citet{Son01} and \citet{Maa05}]:
in the first stage the algorithm searches for small
rectangles with nonzero probability mass; in the second stage those
nonzero probability masses are estimated by maximizing the log
likelihood with a reduced number of unknown quantities}.
\citet{SunWanSun06} and \citet{WuGao11} adopted the same idea
used by \citet{ShiLou95} {to study the joint distribution of
CDF for bivariate interval censored data with Copula models}.

This paper studies bivariate current status data, a special type of
bivariate interval censored data.
Let $(T_1,T_2)$ be the two event times of interest and $(C_1, C_2)$
the two corresponding random monitoring times. In this setting, the
observation of bivariate current status data consists of
%
\begin{equation}
\label{random_biva} X=\bigl(C_1,C_2,
\Delta_1=I(T_1\leq C_1),
\Delta_2=I(T_2\leq C_2)\bigr),
\end{equation}
where $I(\cdot)$ is the indicator function. For bivariate current
status data, \citet{WanDin00} adopted the same approach proposed
by \citet{ShiLou95} to study the association between {the
onset} times of hypertension and diabetes for {Taiwanese} in {a
demographic screening study}. In a~study on HIV transmission,
\citet{JewvanLei05} investigated the relationship
between the time to HIV infection to {the other} partner and the time to
diagnosis of AIDS for the index case by {studying some} smooth
functionals of the marginal CDFs. In both examples, the bivariate
event times have the same monitoring time, that is, $C_1=C_2=C$.
In this paper, we propose a tensor spline-based sieve maximum
likelihood estimation of the joint CDF for bivariate current status
data in a general scenario in which $C_1$ and $C_2$ are allowed to
be different. The proposed method is shown to have a rate of
convergence {potentially} better than $n^{1/3}$ and it can
simultaneously estimate the two marginal CDFs along with the joint
CDF.

The rest of the paper is organized as follows. Section \ref{sec2}
characterizes the spline-based sieve MLE
$\hat{\tau}_n=(\hat{F}_n,\hat{F}_{n,1},\hat{F}_{n,2})$, where
$\hat{F}_n$ is the tensor spline-based estimator of the joint CDF,
and $\hat{F}_{n,1}$ and $\hat{F}_{n,2}$ are the spline-based
estimators of the two corresponding marginal CDFs. Section \ref{sec3}
presents two asymptotic properties (consistency and convergence
rate) of the proposed spline-based sieve MLE. Section \ref{sec4} discusses
the computation of the spline-based estimators. Section \ref{sec5} carries
out a set of simulation studies to examine the finite sample
performance of the proposed method and compares it to the
{conventional NMPLE computed with the} algorithm proposed by
\citet{Maa05}. Section \ref{sec6} summarizes our findings and discusses
some related problems. Section \ref{sec7} provides proofs of the theorems
stated in the early section. Details of some technical lemmas that
are used for proving the theorem and their proofs are included in a
supplementary file.

\section{Tensor spline-based sieve maximum likelihood estimation method}\label{sec2}

\subsection{Maximum likelihood estimation}

{We} consider a sample of $n$ i.i.d. bivariate current status data
denoted in (\ref{random_biva}), $\{(c_{1,k}, c_{2,k}, \delta_{1,k},
\delta_{2,k})\dvtx  k=1,2,\ldots,\allowbreak n\}$. Suppose that $ (T_{1}, T_{2})$
and $ (C_{1}, C_{2})$ are independent. Then the log-likeli\-hood for
the observed data can be expressed by
%
\begin{eqnarray}
\label{like} %
l_n(\cdot;\mbox{data})&=&\sum
_{k=1}^{n} \bigl\{\delta_{1,k}
\delta_{2,k}\log P(T_1 \leq c_{1,k},
T_2 \leq c_{2,k})
\nonumber\\
&&\hphantom{\sum
_{k=1}^{n} \bigl\{}{}+ \delta_{1,k}(1-\delta_{2,k})\log P(T_1 \leq
c_{1,k}, T_2 > c_{2,k})
\nonumber
\\[-8pt]
\\[-8pt]
\nonumber
&&\hphantom{\sum
_{k=1}^{n} \bigl\{}{}+(1-\delta_{1,k})\delta_{2,k}\log P(T_1 >
c_{1,k}, T_2 \leq c_{2,k})
\\
&&\hphantom{\sum
_{k=1}^{n} \bigl\{}{}+(1-\delta_{1,k}) (1-\delta_{2,k})\log P(T_1 >
c_{1,k}, T_2 > c_{2,k})\bigr\}. \nonumber
\end{eqnarray}

Denote $F$ the joint CDF of event times $(T_1,T_2)$ and $F_{1}$ and
$F_{2}$ the marginal CDFs of $F$, respectively. The log-likelihood
(\ref{like}) can be rewritten as
\begin{eqnarray}
\label{flike} %
l_n(F,F_{1},F_{2};
\mbox{data})&=&\sum_{k=1}^{n} \bigl\{
\delta_{1,k}\delta_{2,k}\log F(c_{1,k},
c_{2,k})
\nonumber\\
&&\hphantom{\sum_{k=1}^{n} \bigl\{}{}+ \delta_{1,k}(1-\delta_{2,k})\log\bigl(F_{1}(c_{1,k})-F(c_{1,k},
c_{2,k})\bigr)
\\
&&\hphantom{\sum_{k=1}^{n} \bigl\{}{}+(1-\delta_{1,k})\delta_{2,k}\log\bigl(F_{2}(c_{2,k})-F(c_{1,k},
c_{2,k})\bigr)
\nonumber
\\
&&\hphantom{\sum_{k=1}^{n} \bigl\{}{}+(1-\delta_{1,k}) (1-\delta_{2,k})\log \bigl(1-F_{1}(c_{1,k})-F_{2}(c_{2,k})
\nonumber\\
&&\hspace*{154pt}\hphantom{\sum_{k=1}^{n} \bigl\{}{}+F(c_{1,k},
c_{2,k})\bigr)\bigr\}.\nonumber
\end{eqnarray}

A class of real-valued functions defined in a bounded region
$[L_1,U_1]\times[L_2,U_2]$ is denoted by
\[
\mathcal{F}=\bigl\{\bigl(F(s,t),F_1(s),F_2(t)\bigr)
\dvtx \mbox{ for } (s,t) \in [L_1,U_1]
\times[L_2,U_2] \bigr\},
\]
where $F$, $F_1$ and $F_2$ satisfy the following conditions in (\ref{bvdis}):
%
\begin{eqnarray}
\label{bvdis} %
0&\leq& F(s,t),
\nonumber\\
F\bigl(s',t\bigr)&\leq& F\bigl(s'',t
\bigr),
\nonumber\\
F\bigl(s,t'\bigr)&\leq& F\bigl(s,t''
\bigr),
\nonumber\\
\bigl[F\bigl(s'',t''
\bigr)-F\bigl(s',t''\bigr)\bigr]-
\bigl[(F\bigl(s'',t'\bigr)-F
\bigl(s',t'\bigr)\bigr]&\geq&0,
\nonumber\\
F_1(s)-F(s,t)&\geq&0,
\\
F_2(t)-F(s,t)&\geq&0,
\nonumber\\
\bigl[F_{1}\bigl(s''
\bigr)-F_{1}\bigl(s'\bigr)\bigr]-\bigl[F
\bigl(s'',t\bigr)-F\bigl(s',t\bigr)
\bigr]&\geq&0,
\nonumber\\
\bigl[F_{2}\bigl(t''
\bigr)-F_{2}\bigl(t'\bigr)\bigr]-\bigl[F
\bigl(s,t''\bigr)-F\bigl(s,t'\bigr)
\bigr]&\geq&0,
\nonumber\\
\bigl[1-F_{1}(s)\bigr]-\bigl[F_{2}(t)-F(s,t)\bigr]&\geq&0\nonumber
\end{eqnarray}
for $s'\leq s''$ with $s'$ and $s''$ on $[L_1,U_1]$, and $t'\leq t''$
with $t'$ and $t''$ on $[L_2,U_2]$.

It can be easily argued that if $F$ is a joint CDF and $F_1$ and
$F_2$ are its two corresponding marginal CDFs, $(F,F_1,F_2)\in
{\mathcal{F}}$. On the other hand, for any $(F,F_1,F_2)\in
{\mathcal{F}}$ there exists a bivariate distribution such that $F$
is the joint CDF and $F_1$ and $F_2$ are its two marginal CDFs.
Throughout this paper, $F_0,F_{0,1}$ and $F_{0,2}$ are denoted for
the true joint and marginal CDFs, respectively. The NPMLE for
$(F_0,F_{0,1},F_{0,2})$ is defined as
%
\begin{equation}
\label{nonpa} (\hat{F}_n,\hat{F}_{n,1},\hat{F}_{n,2})=\mathop{
\arg\max}_{(F,F_1,F_2)\in
\mathcal{F}}l_n(F,F_{1},F_{2};
\mbox{data}).
\end{equation}

{The NPNLE of (\ref{nonpa}) is, in general, a challenging problem
both numerically and theoretically. The conventional NPMLE of $F$ is
constructed by a larger number of unknown\vadjust{\goodbreak} quantities representing
the masses in small rectangles. Solving for the NPMLE needs to
perform a constrained high-dimensional nonlinear optimization
[\citet{BetFin99}, \citet{WonYu99}, \citet{GenVan01}, \citet{Son01}, \citet{Maa05}]. Though the conventional
NPMLE of (\ref{nonpa}) can be efficiently computed using the
algorithm developed by \citet{Maa05}, it is, however, well known
that the conventional NPMLE is not uniquely determined
[\citet{Son01},
\citet{Maa05}]. In an unpublished Ph.D. dissertation, \citet{Son01}
showed that the conventional NPMLE of joint CDF with bivariate
current status data can achieve a global rate of convergence of
$n^{3/10}$ in Hellinger distance, which is slightly slower than that
of the NPMLE with univariate current status data.}

This paper adopts a popular dimension reduction method {through
spline-based sieve maximum likelihood estimation.} The main idea of
the spline-based sieve method is to solve problem (\ref{nonpa}) in a
subclass of $\mathcal{F}$ that ``approximates'' to $\mathcal{F}$ when
sample size enlarges. The advantages of the proposed method are that
the spline-based sieve MLE is unique, and it is easy to compute and
analyze. The univariate spline-based sieve MLEs for various models
were {studied} by \citet{She98}, Lu, Zhang and Huang (\citeyear{LuZhaHua07,LuZhaHua09}),
\citet{ZhaHuaHua10} and \citet{Lu10}. Other problems related
to applications of univariate shape-constrained spline estimations
have recently been studied as well. For example, \citet{Mey08}
studied the inference using shape-restricted regression spline
functions and \citet{WanShe10} studied $B$-spline approximation
for a monotone univariate regression function based on grouped data.
For analyzing bivariate distributions, the tensor spline approach
[\citet{deB01}] has been studied by \citet{Sto94} in a~nonparametric
regression setting, by \citet{Koo96} and \citet{Sco92} in a multivariate
density estimation without censored data and, as noted in Section~\ref{sec1},
by \citet{Koo98} in the bivariate density estimation with
bivariate right censored data. Recently, an application of the
tensor $B$-spline estimation of a bivariate monotone function has
also been investigated by \citet{WanTay} in a biomedical
study.

In this paper, we propose a partially monotone tensor spline
estimation of the joint CDF. To solve problem (\ref{nonpa}), the
unknown joint CDF is approximated by a linear combination of the
tensor spline basis functions, and its two marginal CDFs are
approximated by linear combinations of spline basis functions as
well. Then the problem converts to maximizing the sieve log
likelihood with respect to the unknown spline coefficients subject
to a set of inequality constraints.

\subsection{$B$-spline-based estimation}
In this section, the spline-based sieve maximum likelihood
estimation problem is reformulated as a constrained optimization
problem with respect to the coefficients of $B$-spline functions.

Suppose two sets of the normalized $B$-spline basis functions of
order~$l$ [\citet{Sch81}], $\{N_i^{(1),l}(s)\}_{i=1}^{p_n}$ and
$\{N_j^{(2),l}(t)\}_{j=1}^{q_n}$ are constructed in
$[L_1, U_1]\times[L_2,U_2]$ with the knot sequence
$\{u_i\}_{i=1}^{p_n+l}$ satisfying
$L_1=u_1=\cdots=u_l<u_{l+1}<\cdots<u_{p_n}<u_{p_n+1}=\cdots=u_{p_n+l}=U_1$
and the knot sequence $\{v_j\}_{j=1}^{q_n+l}$ satisfying
$L_2=v_1=\cdots=v_l<v_{l+1}<\cdots<v_{q_n}<v_{q_n+1}=\cdots=v_{q_n+l}=U_2$,
where $p_n=O(n^v)$ and $q_n=O(n^v)$ for some $0<v<1$.

Define
\begin{eqnarray*}
\Omega_n&=&\Biggl\{\tau_n=(F_n,F_{n,1},F_{n,2})
\dvtx F_n(s,t)=\sum_{i=1}^{p_n}
\sum_{j=1}^{q_n}\alpha_{i,j}
N^{(1),l}_i(s)N^{(2),l}_j(t),
\\
&&\hspace*{44pt}\hphantom{\Biggl\{}{}F_{n,1}(s)=\sum_{i=1}^{p_n}
\beta_{i}N^{(1),l}_i(s), \ \ F_{n,2}(t)=
\sum_{j=1}^{q_n}\gamma_{j}N^{(2),l}_j(t)\Biggr\},
\end{eqnarray*}
with $\underline{\alpha}=(\alpha_{1,1},\ldots,
\alpha_{p_n,q_n})$, $\underline{\beta}=(\beta_1,\ldots,
\beta_{p_n})$ and $\underline{\gamma}=(\gamma_1,
\ldots,\gamma_{q_n})$
subject to the following conditions in (\ref{consb1}):
%
\begin{eqnarray}
\label{consb1} %
&&\alpha_{1,1}\geq 0,
\nonumber\\
&&\alpha_{1,j+1}-\alpha_{1,j}\geq0 \qquad\mbox{for } j=1,
\ldots,q_n-1,
\nonumber\\
&&\alpha_{i+1,1}-\alpha_{i,1}\geq0 \qquad\mbox{for } i=1,
\ldots,p_n-1,
\\
&&(\alpha_{i+1,j+1}-\alpha_{i+1,j})-(\alpha_{i,j+1}-
\alpha_{i,j})\geq 0
\nonumber\\
&&\qquad\mbox{for } i=1,\ldots,p_n-1, j=1,
\ldots,q_n-1,
\nonumber\\
&&\beta_1-\alpha_{1,q_n}\geq0, \qquad\gamma_1-
\alpha_{p_n,1}\geq 0,
\nonumber\\
&&(\beta_{i+1}-\beta_i)-(\alpha_{i+1,q_n}-
\alpha_{i,q_n})\geq0 \qquad\mbox{for } i=1,\ldots,p_n-1,
\nonumber\\
&&(\gamma_{j+1}-\gamma_j)-(\alpha_{p_n,j+1}-
\alpha_{p_n,j})\geq0 \qquad\mbox{for } j=1,\ldots,q_n-1,
\nonumber\\
&&\beta_{p_n}+\gamma_{q_n}-\alpha_{p_n,q_n}\leq1.\nonumber
\end{eqnarray}
(\ref{consb1}) is established corresponding to the constraints given
in (\ref{bvdis}). Using the properties of $B$-spline, a
straightforward algebra yields $\Omega_n\subset\mathcal{F}$. To
obtain the tensor $B$-spline-based sieve likelihood with bivariate
current status data, $\tau_n=(F_n,F_{n,1},F_{n,2})\in\Omega_n$ is
substituted into (\ref{flike}) to result in
\begin{eqnarray}
\label{likeb} %
&&\tilde{l}_n(\underline{
\alpha},\underline{\beta},\underline{\gamma };\mbox{data})\nonumber\\
&&\qquad= \sum
_{k=1}^{n} \Biggl\{\delta_{1,k}
\delta_{2,k}\log\sum_{i=1}^{p_n}\sum
_{j=2}^{q_n}\alpha_{i,j}N_{i}^{(1),l}(c_{1,k})N_{j}^{(2),l}(c_{2,k})
\nonumber\\
&&\hspace*{20pt}\qquad\quad{}+ \delta_{1,k}(1-\delta_{2,k})\log \Biggl\{ \sum
_{i=1}^{p_n}\beta_{i}N_{i}^{(1),l}(c_{1,k})
\nonumber\\
&&\hspace*{125pt}\qquad{}-\sum_{i=1}^{p_n}\sum
_{j=1}^{q_n}\alpha_{i,j}N_{i}^{(1),l}(c_{1,k})N_{j}^{(2),l}(c_{2,k})
\Biggr\}
\nonumber\\
&&\hspace*{20pt}\qquad\quad{} +(1-\delta_{1,k})\delta_{2,k}\log \Biggl\{ \sum
_{j=1}^{q_n}\gamma_{j}N_{j}^{(2),l}(c_{2,k})
\\
&&\hspace*{125pt}\qquad{}-\sum_{i=1}^{p_n}\sum
_{j=1}^{q_n}\alpha_{i,j}N_{i}^{(1),l}(c_{1,k})N_{j}^{(2),l}(c_{2,k})
\Biggr\}
\nonumber\\
&&\hspace*{20pt}\qquad\quad{} +(1-\delta_{1,k}) (1-\delta_{2,k})\log \Biggl\{1-
\sum_{i=1}^{p_n}\beta_{i}N_{i}^{(1),l}(c_{1,k})
\nonumber\\
&&\hspace*{140pt}\qquad\quad{}- \sum_{j=1}^{q_n}
\gamma_{j}N_{j}^{(2),l}(c_{2,k})\nonumber\\
&&\hspace*{153pt}\qquad{} +\sum
_{i=1}^{p_n}\sum_{j=1}^{q_n}
\alpha_{i,j}N_{i}^{(1),l}(c_{1,k})N_{j}^{(2),l}(c_{2,k})
\Biggr\} \Biggr\}. \nonumber
\end{eqnarray}

Hence, the proposed sieve MLE with the $B$-spline basis functions is
the maximizer of (\ref{likeb}) over $\Omega_n$.

\begin{rem}
The spline-based sieve MLE in $\Omega_n$ is the MLE defined in a
sub-class of $\mathcal{F}$. Hence, the spline-based sieve MLE is
anticipated to have good asymptotic properties if this sub-class
``approximates'' to $\mathcal{F}$ as $n\rightarrow\infty$.
\end{rem}

\section{Asymptotic properties}\label{sec3}
In this section, we describe asymptotic properties of the tensor
spline-based sieve MLE of joint CDF with bivariate current status
data. {Study of the} asymptotic properties of the proposed
sieve estimator requires some regularity conditions, regarding the
event times, observation times and the choice of knot sequences. The
following conditions sufficiently guarantee the results in the
forthcoming theorems.\vspace*{1pt}

\textit{Regularity conditions}:
\begin{longlist}[(C1)]
\item[(C1)] Both $\frac{\partial F_0(s,t)}{\partial s}$ and $\frac{\partial
F_0(s,t)}{\partial t}$ have positive lower bounds in $[L_1,U_1]\times
[L_2,U_2]$.\vspace*{1pt}
\item[(C2)]$\frac{\partial^2 F_0(s,t)}{\partial s\,\partial t}$ has a positive
lower bound $b_0$ in $[L_1,U_1]\times[L_2,U_2].$\vspace*{1pt}
\item[(C3)] $F_0(s,t), F_{0,1}(s)$ and $F_{0,2}(t)$ are all continuous
differentiable up to order $p$
in domain $[L_1,U_1]\times[L_2,U_2]$, $[L_1,U_1]$ and $[L_2,U_2]$,
respectively.\vspace*{1pt}
%
\item[(C4)] The observation times $(C_1,C_2)$ follow a bivariate
distribution defined in $[l_1, u_1]\times[l_2, u_2]$, with $l_1>L_1,
u_1<U_1, l_2>L_2 \mbox{ and } u_2<U_2$.\vspace*{1pt}
\item[(C5)] The density of {the joint distribution of $(C_1,C_2)$} has a
positive lower bound in $[l_1, u_1]\times[l_2, u_2]$.\vspace*{1pt}
\item[(C6)]
The knot sequences $\{u_i\}_{i=1}^{p_n+l}$ and
$\{v_j\}_{j=1}^{q_n+l}$ of the $B$-spline basis functions,
$\{N^{(1),l}_i\}_{i=1}^{p_n}$ and $\{N^{(2),l}_j\}_{j=1}^{q_n}$,
satisfy that both
$\frac{\min_{i}\Delta_{i}^{(u)}}{\max_{i}\Delta_{i}^{(u)}}$ and
$\frac{\min_{j}\Delta_{j}^{(v)}}{\max_{j}\Delta_{j}^{(v)}}$ have
positive lower bounds, where $\Delta_i^{(u)}=u_{i+1}-u_i$ for
$i=l,\ldots,p_n$ and
$\Delta_j^{(v)}=v_{j+1}-v_j$ for $j=l,\ldots,q_n$.\vspace*{-1pt}
\end{longlist}
%
\begin{rem}
(C1) implies that $\frac{d F_{0,1}(s)}{d s}$ and $\frac{d F_{0,2}(t)}{d
t}$ have positive lower bounds on $[L_1,U_1]$ and $[L_2,U_2]$,
respectively. (C3) implies that both $\frac{\partial F_0(s,t)}{\partial
s}$ and $\frac{\partial F_0(s,t)}{\partial t}$ have positive upper
bounds in $[L_1,U_1]\times[L_2,U_2]$; $\frac{d F_{0,1}(s)}{d s}$ and
$\frac{d F_{0,2}(t)}{d t}$ have positive upper bounds on $[L_1,U_1]$
and $[L_2,U_2]$, respectively.\vspace*{-1pt}
\end{rem}
%

Let
\begin{eqnarray*}
\Omega_{n,1}&=&\Biggl\{\tau=(F_n,F_{n,1},F_{n,2})
\dvtx F_n(s,t)=\sum_{i=1}^{p_n}
\sum_{j=1}^{q_n}\alpha_{i,j}
N^{(1),l}_i(s)N^{(2),l}_j(t),
\\[-2pt]
&&\hspace*{44pt}\hphantom{\Biggl\{}F_{n,1}(s)=\sum_{i=1}^{p_n}
\beta_{i}N^{(1),l}_i(s), F_{n,2}(t)=\sum
_{j=1}^{q_n}\gamma_{j}N^{(2),l}_j(t)\Biggr\},\vspace*{1pt}
\end{eqnarray*}
with $\underline{\alpha}=(\alpha_{1,1},\ldots,
\alpha_{p_n,q_n})$, $\underline{\beta}=(\beta_1,\ldots,
\beta_{p_n})$ and $\underline{\gamma}=(\gamma_1,
\ldots,\gamma_{q_n})$
subject to the following conditions in (\ref{consb2}):
%
\begin{eqnarray}
\label{consb2} %
&&\alpha_{1,1}\geq0,
\nonumber\\
&&\alpha_{1,j+1}-\alpha_{1,j}\geq0 \qquad\mbox{for } j=1,
\ldots,q_n-1,
\nonumber\\
&&\alpha_{i+1,1}-\alpha_{i,1}\geq0 \qquad\mbox{for } i=1,
\ldots,p_n-1,
\nonumber\\
&&(\alpha_{i+1,j+1}-\alpha_{i+1,j})-(\alpha_{i,j+1}-
\alpha_{i,j})
\nonumber\\
&&\qquad\geq \frac{b_0\min_{i_1: l\leq i_1\leq p_n }\Delta_{i_1}^{(u)}\min_{j_1:
l\leq j_1\leq q_n}\Delta_{j_1}^{(v)}}{l^2}
\nonumber
\\[-9pt]
\\[-9pt]
\nonumber
&&\qquad\quad\mbox{for } i=1,\ldots,p_n-1, j=1,
\ldots,q_n-1,
\nonumber\\
&&\beta_1-\alpha_{1,q_n}\geq0, \qquad\gamma_1-
\alpha_{p_n,1}\geq0,
\nonumber\\
&&(\beta_{i+1}-\beta_i)-(\alpha_{i+1,q_n}-
\alpha_{i,q_n})\geq0\qquad \mbox{for } i=1,\ldots,p_n-1,
\nonumber\\
&&(\gamma_{j+1}-\gamma_j)-(\alpha_{p_n,j+1}-
\alpha_{p_n,j})\geq0 \qquad\mbox{for } j=1,\ldots,q_n-1,
\nonumber\\
&&\beta_{p_n}+\gamma_{q_n}-\alpha_{p_n,q_n}\leq1.\nonumber
\end{eqnarray}


\begin{rem}
Note that $\Omega_{n,1}$ is a sub-class of $\Omega_{n}$ due to the
change from the forth inequality of (\ref{consb1}) to that of
(\ref{consb2}). The choice of $\Omega_{n,1}$ is mainly for the
technical convenience in justifying the asymptotic properties. In
the forth inequality of~(\ref{consb2}), $b_0$ is the positive lower
bound of $\frac{\partial^2 F_0(s,t)}{\partial s\,\partial t}$ stated
in~(C2). This inequality will guarantee that $\frac{\partial^2
F_n(s,t)}{\partial s\,\partial t}$ also has a positive lower bound
which is
necessary for the proof of Lemma 0.1 in the supplemental article [Wu and Zhang (\citeyear{WuZha})].
It is obvious
that as sample $n$ increases to infinity, the right-hand side of the
forth inequality in~(\ref{consb2}) will approach to~0.\vadjust{\goodbreak}
\end{rem}

We study the asymptotic properties in the feasible region of the
observation times: $ [l_1, u_1]\times[l_2, u_2]$. Let
$\Omega_n^{\prime}=\{\tau\in\Omega_{n,1} , \mbox{ for } (s,t)\in[l_1,
u_1]\times[l_2, u_2]\}$ and let
$\tau_0=(F_0(s,t),F_{0,1}(s),F_{0,2}(t))$ with $(s,t)\in[l_1,
u_1]\times[l_2, u_2]$. Under~(C4), the maximization of
$\tilde{l}_n(\underline{\alpha},\underline{\beta},\underline{\gamma
};\mbox{data})$
over $\Omega_{n,1}$ is actually the maximization of
$\tilde{l}_n(\underline{\alpha},\underline{\beta},\underline{\gamma
};\mbox{data})$
over $\Omega_n^{\prime}$. Throughout the study of asymptotic properties,
we denote $\hat{\tau}_n$ the maximizer of
$\tilde{l}_n(\underline{\alpha},\underline{\beta},\underline{\gamma
};\mbox{data})$
over $\Omega_n^{\prime}$.

Denote $L_r(Q)$ the norm associated with a probability measure $Q$
which is defined as
\[
\|f\|_{L_r(Q)}=\bigl(Q|f|^r\bigr)^{1/r}= \biggl(
\int|f|^r\,dQ \biggr)^{1/r}.
\]
In the following, $L_r(P_{C_1,C_2})$, $L_r(P_{C_1})$ and
$L_r(P_{C_2})$ are denoted as the $L_r$-norms associated with the
joint and marginal probability measures of the observation times
$(C_1,C_2)$, respectively, and $L_r(P)$ is denoted as the
$L_r$-norm associated with the joint probability measure $P$ of
observation and event times $(T_1,T_2,C_1,C_2)$.

Based on the $L_2$-norms, the distance between $\tau_n=(F_n,F_{n,1},F_{n,2})\in\Omega_n^{\prime}$ and $\tau_0=(F_0,F_{0,1},F_{0,2})$ is defined as
\begin{eqnarray*}
&& d(\tau_n,\tau_0)\\
&&\qquad=\bigl(\|F_n-F_0
\|^2_{L_2(P_{C_1,C_2})} +\|F_{n,1}-F_{0,1}
\|^2_{L_2(P_{C_1})} +\|F_{n,2}-F_{0,2}
\|^2_{L_2(P_{C_2})}\bigr)^{1/2}.\vspace*{-1pt}
\end{eqnarray*}


\begin{thmm}\label{th3.1}
Suppose \textup{(C2)--(C6)} hold, and $p_n=O(n^v)$, $q_n=O(n^v)$ for $v<1$; that
is, the numbers of interior knots of knot sequences
$\{u_i\}_1^{p_n+l}$ and $\{v_j\}_1^{q_n+l}$ are both in the order of
$n^v$ for $v<1$. Then
\[
d(\hat{\tau}_n,\tau_0)\rightarrow_p0
\qquad\mbox{as } n\rightarrow\infty.\vspace*{-1pt}
\]
\end{thmm}

\begin{thmm}\label{th3.2}
Suppose \textup{(C1)--(C6)} hold, and $p_n=O(n^v)$, $q_n=O(n^v)$ for $v< 1$; that
is, the numbers of interior knots of knot sequences
$\{u_i\}_1^{p_n+l}$ and $\{v_j\}_1^{q_n+l}$ are both in the order of
$n^v$ for $v< 1$. Then
\[
d(\hat{\tau}_n,\tau_0)=O_p
\bigl(n^{-\min\{pv,(1-2v)/2\}}\bigr).\vspace*{-1pt}
\]
\end{thmm}

\begin{rem}
Theorem \ref{th3.2} {implies that the optimal rate of convergence of the
proposed estimator is $n^{{p}/{(2(p+1))}}$, achieved by letting
$pv=(1-2v)/2$. This rate is equal to $n^{1/3}$ when $p=2$ and
improves as $p$ (the degree of smoothness of the true joint
distribution) increases. Nonetheless, the rate will never exceed
$n^{1/2}$. The result of Theorem \ref{th3.2} also indicates that the
proposed method potentially results in an estimate of the targeted
joint CDF with a faster convergence rate than the conventional NPMLE
method given in \citet{Son01}}.\vspace*{-1pt}
\end{rem}

\section{Computation of the spline-based sieve MLE}\label{sec4}
{For the $B$-spline-based sieve MLE, the constraint set
(\ref{consb2})\vadjust{\goodbreak} complicates the numerical implementation}. We propose
to compute the sieve MLE using $I$-spline basis functions for the
sake of numerical convenience. The $I$-spline basis functions are
defined by \citet{Ram88} as
%
\begin{equation}
\label{ramsay}\qquad I_i^{l}(s)= %
\cases{ 0, &\quad $i>j,$
\vspace*{2pt}
\cr
\displaystyle\sum_{m=i}^j(u_{m+l+1}-u_m)M_m^{l+1}(s)/(l+1),&\quad
$j-l+1\leq i\leq j,$\vspace*{2pt}
\cr
1,& \quad $i<j-l+1$ } %
\end{equation}
for $u_j\leq s<u_{j+1}$, where $M_m^{l}$s are the $M$-spline basis
functions of order $l$, studied by
\citet{CurSch66}, and can be calculated recursively by\vspace*{1pt}
\begin{eqnarray*}
\label{ms1}
M_i^{1}(s)&=& \frac{1}{u_{i+1}-u_i},\qquad u_i\leq s<u_{i+1},
\\
M_i^l(s)&=&\frac
{l[(s-u_i)M_i^{l-1}(s)+(u_{i+l}-s)M_{i+1}^{l-1}(s)]}{(l-1)(u_{i+l}-u_i)}.
\end{eqnarray*}
By the relationship between the $B$-spline basis functions and the
$M$-spline basis functions [\citet{Sch81}], it can be easily
{demonstrated} that the $I$-spline basis function defined
{in} (\ref{ramsay}) can be expressed by a sum of the
$B$-spline basis functions\vspace*{1pt}
%
\begin{equation}
\label{3btoi1} I_i^{l-1}(s)=\sum
_{m=i}^{p_n} N_m^{l}(s).
\end{equation}
{Consequently}, the spline-based sieve {space} can be
re{constructed using} the $I$-spline basis functions
{with a different set of constraints:}\vspace*{1pt}
\begin{eqnarray*}
\Theta_n&=&\Biggl\{\tau_n=(F_n,F_{n,1},F_{n,2})
\dvtx F_n(s,t)=\sum_{i=1}^{p_n}
\sum_{j=1}^{q_n}\eta_{i,j}
I^{(1),l-1}_i(s)I^{(2),l-1}_j(t),
\\
&&\hphantom{\Biggl\{}F_{n,1}(s)=\sum_{i=1}^{p_n}\Biggl
\{\sum_{j=1}^{q_n}\eta_{i,j}+
\omega_{i}\Biggr\} I_{i}^{(1),l-1}(s),
\\
&&\hspace*{111pt}\hphantom{\Biggl\{}F_{n,2}(t)=\sum_{j=1}^{q_n}\Biggl
\{\sum_{i=1}^{p_n}\eta_{i,j}+
\pi_{j}\Biggr\} I_{j}^{(2),l-1}(t)\Biggr\}
\end{eqnarray*}
with $\underline{\eta}=(\eta_{1,1},\ldots,\eta_{p_n,q_n})$,
$\underline{\omega}=(\omega_1,\ldots,\omega_{p_n})$ and
$\underline{\pi}=(\pi_1,\ldots,\pi_{q_n})$
subject to the following conditions in (\ref{consi}),\vspace*{1pt}
\begin{eqnarray}
\label{consi} %
&&\eta_{i,j}\geq0 \qquad\mbox{for } i=1,
\ldots, p_n , j=1,\ldots,q_n,
\nonumber\\
&&\omega_i\geq0,\qquad i=1,\ldots,p_n,
\\
&&\pi_j\geq0,\qquad j=1,\ldots,q_n,
\nonumber\\
&&\sum_{i=1}^{p_n}\sum
_{j=1}^{q_n}\eta_{i,j}+\sum
_{i=1}^{p_n}\omega_i+ \sum
_{j=1}^{q_n}\pi_j\leq1. \nonumber%
\end{eqnarray}
Then the spline-based sieve log likelihood can be also expressed in
$I$-spline, and the spline-based sieve MLE can be obtained by
maximizing the log likelihood in $I$-spline over $\Theta_n$.

\begin{rem}
Class $\Theta_n$ is actually equivalent to $\Omega_n$, and hence the
$I$-spline-based sieve MLE is the same as the $B$-spline-based sieve
MLE. It is advocated in numerical implementation due to the
simplicity of the constraints in class $\Theta_n$.
\end{rem}



Given $p_n$ and $q_n$, the proposed sieve estimation problem
described above is actually a restricted parametric maximum
likelihood estimation problem with respect to the coefficients
{associated with the} $I$-spline and the tensor $I$-spline basis
functions. \citet{Jam04} generalized the gradient projection
algorithm originally proposed by \citet{Ros60} using a weighted
$L_2$-norm $\|x\|=x'Wx$ with a positive definite matrix $W$ for the
restricted maximum likelihood estimation problems.
{Because the constraint set (\ref{consi}) is made by linear
inequalities, the maximization of~(\ref{flike}) in the $I$-spline
form over $\Theta_n$ can be efficiently implemented by the
generalized gradient projection algorithm [\citet{Jam04}] and is
described as follows}.

First we rewrite (\ref{consi}) as $X\theta\leq y$, where\vspace*{1pt}
$X=(x_1,x_2,\ldots,x_{p_n\cdot q_n+p_n+q_n},\break  x_{p_n\cdot
q_n+p_n+q_n+1})^T$ with $x_1=(-1,0,\ldots,0)^T$,
$x_2=(0,-1,0,\ldots,0)^T$,\break $x_{p_n\cdot
q_n+p_n+q_n}=(0,\ldots,0,-1)^T$, $x_{p_n\cdot
q_n+p_n+q_n+1}=(1,\ldots,1)^T$;
$\theta=(\underline{\eta},\underline{\omega},\underline{\pi})=(\theta_1,\theta_2,\ldots,\theta_{p_n\cdot
q_n+p_n+q_n})$; and $y=(0,\ldots,0,1)^T$. If some $I$-spline
coefficients equal 0 or all coefficients sum up to 1, then we say
their {corresponding} constraints are active and let
$\bar{X}\theta=\bar{y}$ represent all {the} active
constraints {and a vector $\Lambda$ of integers to index the
active constraints}.
For example, if $\Lambda=(2,1,p_n\cdot q_n+p_n+q_n+1)$, then
{the second, first and last constraints become active}, and
$\bar{X}=(x_2,x_1,x_{p_n\cdot q_n+p_n+q_n+1})^T$ {and
$\bar{y}=(0,0,1)^T$}.

Let $\dot{\tilde{l}}(\theta)$ and $H(\theta)$ be the gradient and
Hessian matrix of the log likelihood given by (\ref{flike}) in the
$I$-spline form, respectively. Note that $H(\theta)$ may not be
negative definite for every $\theta$. We use $W=-H(\theta)+\delta
I$, where $I$ is identity matrix, and $\delta>0$ is chosen as any
value that guarantees $W$ being positive definite. With that
introduced, the generalized gradient projection algorithm is
implemented as follows.

\begin{enumerate}[Step 1]
\item[Step 1]
(\textit{Computing the feasible search direction}). Compute
\begin{eqnarray*}
\underline{d}&=&(d_1,d_2,\ldots,d_{p_n\cdot q_n+p_n+q_n})
\\
&=&\bigl\{I-W^{-1}\bar{X}^{T}\bigl(\bar{X}W^{-1}
\bar{X}^{T}\bigr)^{-1}\bar{X}\bigr\} W^{-1}\dot{
\tilde{l}}(\theta).
\end{eqnarray*}
\item[Step 2]
(\textit{Forcing the updated $\theta$ to fulfill the constraints}). Compute
\begin{eqnarray*}
\gamma= %
\cases{\displaystyle\min\biggl\{\min_{i:d_i<0}\biggl\{-
\frac{\theta_i}{d_i}\biggr\},\frac{1-\sum_{i=1}^{p_n\cdot q_n+p_n+q_n}\theta_i}{\sum_{i=1}^{p_n\cdot
q_n+p_n+q_n}\,d_i}\biggr\},\vspace*{2pt}\cr
\hspace*{61pt}\qquad\mbox{if }
\displaystyle\sum_{i=1}^{p_n\cdot q_n+p_n+q_n}\,d_i>0,
\vspace*{2pt}
\cr
\displaystyle\min_{i:d_i<0}\biggl\{-\frac{\theta_i}{d_i}\biggr\},\qquad
\mbox{else}.} %
\end{eqnarray*}
Doing so guarantees that $\theta_i+\gamma d_i\geq0$ for
$i=1,2,\ldots, p_n\cdot q_n+p_n+q_n$, and $\sum_{i=1}^{p_n\cdot
q_n+p_n+q_n}(\theta_i+\gamma d_i)\leq1$.
\item[Step 3]
(\textit{Updating the solution by step-halving line search}). Find the
smallest integer $k$ starting from 0 such that
\[
\tilde{l}_n\bigl(\theta+(1/2)^k\gamma\underline{d};\cdot
\bigr)\geq\tilde {l}_n(\theta;\cdot).
\]
Replace $\theta$ by $\tilde{\theta}=\theta+\min\{(1/2)^k\gamma,0.5\}
\underline{d}$.
\item[Step 4]
(\textit{Updating $\Lambda$ and $\bar{X}$}). Modify $\Lambda$ by adding
indexes of new $I$-spline coefficients when these new coefficients
become 0 and adding $p_n\cdot q_n+p_n+q_n+1$ when the sum of all
$I$-spline coefficients becomes 1. Modify $\bar{X}$ accordingly.
\item[Step 5]
(\textit{Checking the stopping criterion}). If $\|\underline{d}\|\geq
\varepsilon$, for small $\varepsilon$, go to Step 1; otherwise, compute
$\lambda=(\bar{X}W^{-1}\bar{X}^{T})^{-1}\bar{X}W^{-1}\dot{\tilde
{l}}(\theta)$.
\begin{enumerate}[(ii)]
\item[(i)]
If the $j$th component $\lambda_j\geq0$ for all $j$, set $\hat{\theta
}=\theta$ and stop.
\item[(ii)]
If there is at least one $j$ such that $\lambda_j<0$, let $j^{*}=
\arg\min_{j:\lambda_j<0}\{\lambda_j\}$, then remove $j^{*}$th component
from $\Lambda$ and remove the $j^{*}$th row from $\bar{X}$, and go to
Step 1.
\end{enumerate}
\end{enumerate}

\section{Simulation studies}\label{sec5}
Copula models are often used in studying bivariate event time data
[\citet{ShiLou95}, \citet{WanDin00},
\citet{SunWanSun06},
\citet{Zhaetal10}]

We consider the bivariate Clayton copula function
\[
C_{\alpha}(u,v)=\bigl(u^{(1-\alpha)}+v^{(1-\alpha)}-1
\bigr)^{{1}/{(1-\alpha)}},
\]
with $\alpha>1$.
For the Clayton copula, a larger $\alpha$ corresponds to a stronger
positive association between the two random variables. The
association parameter $\alpha$ and Kendall's $\tau$ for the Clayton
copula is related by $\tau=\frac{\alpha-1}{\alpha+1}$.

In the simulation studies,
we compare the proposed sieve MLE to the {conventional} NPMLE,
{computed using} the algorithm {developed} by \citet{Maa05}. As we
mentioned previously, this NPMLE is not unique. Only the total mass
in each selected rectangle is estimated, therefore the estimated
joint CDF is based on where the mass is placed in each rectangle. We
denote U-NPMLE and L-NPMLE as the NPMLE {for which the probability
mass is placed at the upper right and lower left corners of each
rectangle, respectively}.

The proposed sieve MLE and both U-NPMLE and L-NPMLE are evaluated with
various combinations of Kendall's $\tau$ $(\tau=0.25, 0.75)$ and sample
sizes ($n=100, 200$).
Under each of the four settings, the Monte-Carlo simulation with 500
repetitions is conducted, and the cubic ($l=4$) $I$-spline basis
functions are used
in the proposed sieve estimation method.
The event times ($T_1$, $T_2$), monitoring times ($C_1$, $C_2$) and
the knots selection of the cubic $I$-spline basis functions are
{illustrated} as follows:
\begin{longlist}[(iii)]
\item[(i)]
(\textit{Event times}). $(T_1,T_2)$ are generated from the Clayton copula
with the two marginal distributions being exponential with the rate
parameter~0.5.
Under this setting, $\operatorname{Pr}(T_i\geq5)<0.1$ for $i=1,2$ and
$[L_1,U_1]\times[L_2,U_2]$ is chosen to be $[0,5] \times[0,5]$.
\item[(ii)]
(\textit{Censoring times}). Both $C_1$ and $C_2$ are generated
independently from the uniform distribution
on $[0.0201,4.7698]$ [$\operatorname{Pr}(0<T_i<0.0201)=\operatorname{Pr}(4.7698<T_i<5)=0.01$, for
$i=1,2$]. The observation region $[l_1,u_1]\times
[l_2,u_2]=[0.0201,4.7698]\times[0.0201,4.7698]$
is inside $[0,5]\times[0,5]$ and the CDFs are bounded away from 0 and
1 inside the observation region.
\item[(iii)]
(\textit{Knots selection}). As in other spline-based estimations
[Lu, Zhang and Huang (\citeyear{LuZhaHua07}, \citeyear{LuZhaHua09}), \citet{ZhaHuaHua10}
{and \citet{WuGao11}}], the number of interior knots $m_n$ is chosen as
$ [n^{1/3} ]-1$, {where $ [n^{1/3} ]$ is the integer
part of $n^{1/3}$}. For
moderate sample sizes, say $n=100,200$,
our experiments show that $m_n= [n^{1/3} ]-1$ is a
reasonable choice for the number of interior knots. Therefore, we
choose 4 and 5 as the numbers of interior knots for sample sizes 100
and 200, respectively. The number of spline basis functions is
determined by $p_n=q_n=m_n+4$ in our computation. Two end knots of
all knot sequences are chosen to be 0 and 5. For each sample of
bivariate observation times $(C_1,C_2)$, the interior knots for
$\{I_i^{(1),3}\}_{i=1}^{p_n}$ and $\{I_j^{(2),3}\}_{j=1}^{q_n}$ are
allocated at the $k/(m_n+1)$ quantiles ($k=1,\ldots,m_n$) of the
samples of $C_1$ and $C_2$, respectively.
\end{longlist}


\begin{table}[t!]
\caption{Comparison of the estimation bias and
square root of mean square error among the sieve MLE, U-NPMLE and
L-NPMLE at four selected points}\label{tab1}
\begin{tabular*}{\textwidth}{@{\extracolsep{\fill}}lccccccc@{}}
\hline
&&\multicolumn{6}{c@{}}{$\bolds{T_2}$}\\[-6pt]
&&\multicolumn{6}{c@{}}{\hrulefill}\\
&&\multicolumn{3}{c}{\textbf{0.1}}&\multicolumn{3}{c@{}}{\textbf{4.6}}\\[-4pt]
&&\multicolumn{3}{c}{\hrulefill}&\multicolumn{3}{c@{}}{\hrulefill}\\
\multicolumn{1}{@{}l}{$\bolds{T_1}$} && \textbf{Sieve} &\textbf{U-Non}&\textbf{L-Non}& \textbf{Sieve}&\textbf{U-Non}&\textbf{L-Non} \\ %
\hline
\multicolumn{8}{@{}l}{Sample size $n=100$, Kendall's $\tau=0.25$}\\
0.1&$\mathrm{Bias}$& $-$5.00e--3& $-$1.78e--2& $-$1.91e--2& \phantom{$-$}2.69e--2& $-$2.22e--3& $-$3.93e--2\\
&$\mathrm{MSE}^{1/2}$& \phantom{$-$}2.75e--2& \phantom{$-$}2.24e--2& \phantom{$-$}1.91e--2& \phantom{$-$}7.32e--2& \phantom{$-$}6.68e--2& \phantom{$-$}5.53e--2\\[3pt]
4.6&$\mathrm{Bias}$& \phantom{$-$}2.33e--2& $-$2.69e--2& $-$4.19e--2& \phantom{$-$}4.04e--2& \phantom{$-$}1.32e--1& \phantom{$-$}1.09e--1\\
&$\mathrm{MSE}^{1/2}$& \phantom{$-$}7.18e--2& \phantom{$-$}6.68e--2& \phantom{$-$}5.17e--2& \phantom{$-$}8.24e--2& \phantom{$-$}1.49e--1& \phantom{$-$}1.35e--1\\[6pt]
\multicolumn{8}{@{}l}{Sample size $n=200$, Kendall's $\tau=0.25$}\\
0.1&$\mathrm{Bias}$& $-$4.39e--3& $-$1.85e--2& $-$1.91e--2& \phantom{$-$}2.26e--2& $-$2.87e--2& $-$3.89e--2\\
&$\mathrm{MSE}^{1/2}$& \phantom{$-$}2.42e--2& \phantom{$-$}1.98e--2& \phantom{$-$}1.91e--2& \phantom{$-$}6.05e--2& \phantom{$-$}5.52e--2& \phantom{$-$}5.04e--2\\[3pt]
{4.6}&$\mathrm{Bias}$& \phantom{$-$}1.65e--2& $-$3.29e--2& $-$4.03e--2& \phantom{$-$}2.15e--2& \phantom{$-$}1.10e--1& \phantom{$-$}9.65e--2\\
&$\mathrm{MSE}^{1/2}$& \phantom{$-$}5.31e--2& \phantom{$-$}5.50e--2& \phantom{$-$}5.30e--2& \phantom{$-$}6.10e--2& \phantom{$-$}1.29e--1& \phantom{$-$}1.21e--1\\[6pt]
\multicolumn{8}{@{}l}{Sample size $n=100$, Kendall's $\tau=0.75$}\\
{0.1}&$\mathrm{Bias}$& $-$1.81e--2& $-$3.62e--2& $-$4.33e--2& \phantom{$-$}2.91e--2& $-$1.95e--2& $-$4.11e--2\\
&$\mathrm{MSE}^{1/2}$& \phantom{$-$}4.63e--2& \phantom{$-$}5.36e--2& \phantom{$-$}4.34e--2& \phantom{$-$}7.88e--2& \phantom{$-$}8.19e--2& \phantom{$-$}5.62e--2\\[3pt]
{4.6}&$\mathrm{Bias}$& \phantom{$-$}3.08e--2& $-$1.90e--2& $-$4.04e--2& \phantom{$-$}1.98e--2& \phantom{$-$}1.03e--1& \phantom{$-$}8.09e--2\\
&$\mathrm{MSE}^{1/2}$& \phantom{$-$}8.16e--2& \phantom{$-$}8.45e--2& \phantom{$-$}5.83e--2& \phantom{$-$}6.47e--2& \phantom{$-$}1.22e--1&\phantom{$-$}1.08e--1\\[6pt]
\multicolumn{8}{@{}l}{Sample size $n=200$, Kendall's $\tau=0.75$}\\
{0.1}&$\mathrm{Bias}$& $-$2.03e--2& $-$4.00e--2& $-$4.31e--2& \phantom{$-$}2.01e--2& $-$2.48e--2& $-$3.81e--2\\
&$\mathrm{MSE}^{1/2}$& \phantom{$-$}3.86e--2& \phantom{$-$}4.52e--2& \phantom{$-$}4.37e--2& \phantom{$-$}5.87e--2& \phantom{$-$}5.90e--2& \phantom{$-$}5.27e--2\\[3pt]
{4.6}&$\mathrm{Bias}$& \phantom{$-$}2.09e--2& $-$2.48e--2& $-$3.93e--2& \phantom{$-$}1.08e--2& \phantom{$-$}8.37e--2& \phantom{$-$}7.20e--2\\
&$\mathrm{MSE}^{1/2}$& \phantom{$-$}6.00e--2& \phantom{$-$}6.27e--2& \phantom{$-$}5.35e--2& \phantom{$-$}5.24e--2& \phantom{$-$}1.04e--1& \phantom{$-$}9.40e--2\\
\hline
\end{tabular*}
\end{table}

\begin{table}[t!]
\caption{Comparison of the overall estimation
biases and the overall mean square errors among sieve MLE, U-NPMLE
and L-NPMLE}\label{tab2}
\begin{tabular*}{\textwidth}{@{\extracolsep{\fill}}lccccccc@{}}
\hline
&&\multicolumn{6}{c@{}}{\textbf{Sample size}}\\[-4pt]
&&\multicolumn{6}{c@{}}{\hrulefill}\\
&&\multicolumn{3}{c}{\textbf{100}} &\multicolumn{3}{c@{}}{\textbf{200}}\\[-4pt]
&&\multicolumn{3}{c}{\hrulefill} &\multicolumn{3}{c@{}}{\hrulefill}\\
{$\bolds{\tau}$} & & \textbf{Sieve} &\textbf{U-Non}& \textbf{L-Non}&\textbf{Sieve}&\textbf{U-Non}& \textbf{L-Non}\\
\hline
0.25 &$|\mathrm{Bias}|$&7.56e--3&3.24e--2&4.24e--2& 6.70e--3& 2.62e--2&3.11e--2\\
&$\mathrm{MSE}^{1/2}$& 7.93e--2& 1.25e--1&1.26e--1& 6.13e--2&
1.03e--1&1.03e--1\\[3pt]
0.75 &$|\mathrm{Bias}|$& 1.11e--2& 1.50e--2&2.49e--2& 7.29e--3& 1.33e--2&1.88e--2\\
&$\mathrm{MSE}^{1/2}$& 7.40e--2& 1.10e--1&1.11e--1& 5.77e--2& 8.45e--2&8.53e--2\\
\hline
\end{tabular*}
\end{table}

\begin{figure}

\includegraphics{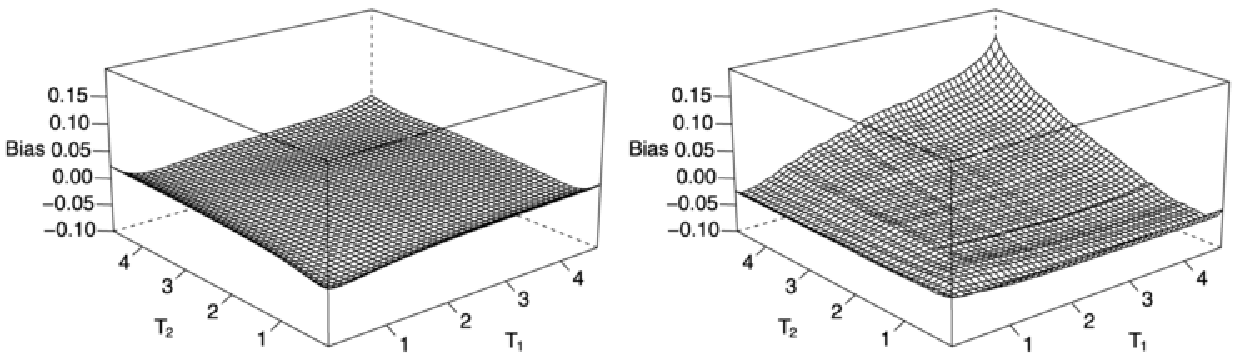}

\caption{Comparison of the estimation bias between the sieve
MLE (left) and the U-NPMLE (right) for the joint CDF when sample size
$n=200$, Kendall's $\tau=0.25$.}\label{fig1}
\end{figure}

\begin{figure}

\includegraphics{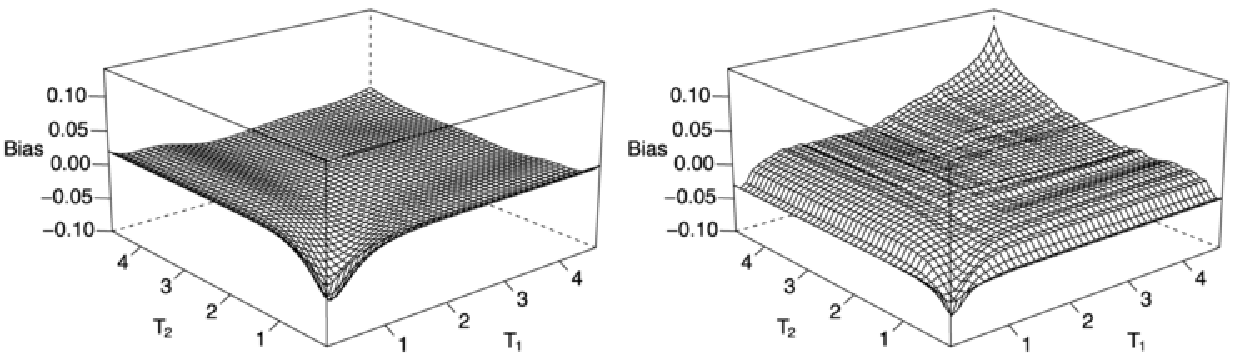}

\caption{Comparison of the estimation bias between the sieve
MLE (left) and the U-NPMLE (right) for the joint CDF when sample size
$n=200$, Kendall's $\tau=0.75$.}\label{fig2}
\end{figure}

\begin{figure}

\includegraphics{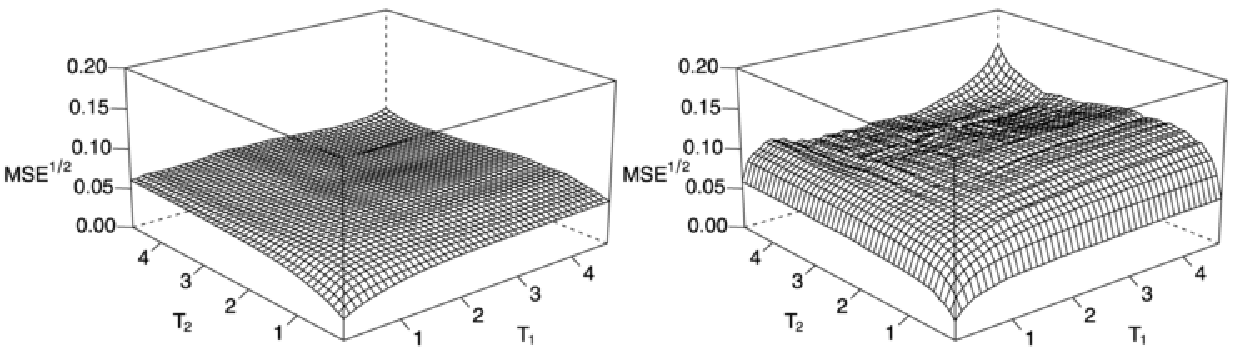}

\caption{Comparison of the square root of mean square error
between the sieve MLE (left) and the U-NPMLE (right) for the joint CDF
when sample size $n=200$, Kendall's $\tau=0.25$.}\label{fig3}
\end{figure}

Table \ref{tab1} displays the estimation biases (\textit{Bias}) and the square
roots of mean square errors ($\mathrm{MSE}^{1/2}$) from the Monte-Carlo
simulation of 500 repetitions for the proposed sieve MLE (\textit{Sieve}) and both U-NPMLE (\textit{U-Non}) and L-NPMLE (\textit{L-Non}) of
the bivariate CDF at 4 {selected} pairs of time points $(s_1,s_2)$
near the corners of the estimation region with different sample
sizes and Kendall's $\tau$ values. The estimation results at those
selected points are comparable among the three estimators. Table \ref{tab2}
{presents} the {overall} estimation bias and mean square error {for
the} three estimators {by calculating the average of absolute values
of estimation bias and the average of square roots of mean square
error taking from} 2209 {pairs} of $(s_1,s_2)$ with both $s_1$ and
$s_2$ {ranging uniformly} from 0.1 to 4.7. It appears that the sieve
MLE outperforms its counterparts {with a} smaller overall bias and a
smaller overall mean square error.
The mean square error of the proposed sieve MLE noticeably decreases
as sample size increases from 100 to 200.

\begin{figure}

\includegraphics{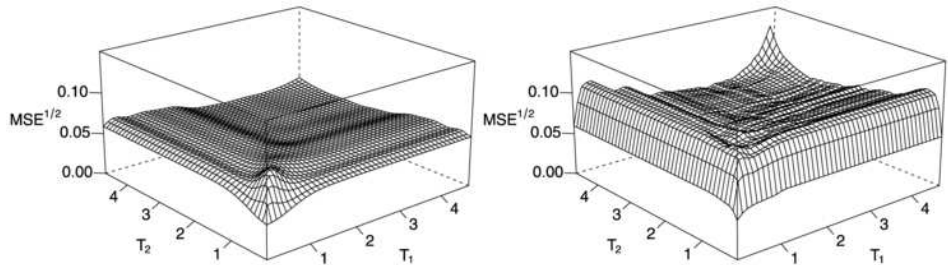}

\caption{Comparison of the square root of mean square error
between the sieve MLE (left) and the U-NPMLE (right) for the joint CDF
when sample size $n=200$, Kendall's $\tau=0.75$.}\label{fig4}
\end{figure}

\begin{figure}

\includegraphics{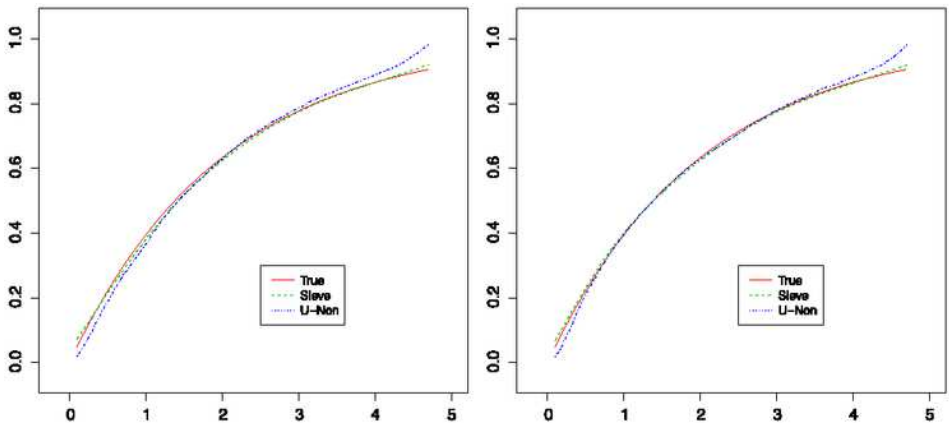}

\caption{Comparison of the estimation bias between the sieve
MLE and the U-NPMLE for estimating the marginal CDF of $T_1$ when
sample size $n=200$
(left: Kendall's $\tau=0.25$; right: Kendall's $\tau=0.75$).}\label{fig5}
\end{figure}

For sample size $n=200$, the estimation biases and the square roots
of mean square error of the sieve MLE and U-NPMLE for the joint CDF
from the same Monte-Carlo simulation are graphed in Figure \ref{fig1} through
Figure~\ref{fig4} for Kendall's $\tau=0.25$ and 0.75. These figures clearly
indicate that the bias and the MSE of the sieve MLE are noticeably
smaller than that of U-NPMLE inside the closed region $[0.1,
4.7]\times[0.1,4.7]$. It is also seen that the bias of the sieve
MLE near the origin increases as Kendall's $\tau$ increases. As
a~by-product of the estimation methods, the {average} estimate of the
marginal CDF of $T_1$ from the same Monte-Carlo simulation for both
the proposed sieve MLE (\textit{Sieve}) and U-NPMLE (\textit{U-Non}) are
also computed and plotted in Figure \ref{fig5} along with the true marginal
CDF (\textit{True}), $F_1$. Figure \ref{fig5} clearly indicates that the bias
of the proposed sieve MLE for the marginal CDF is {markedly} smaller
than that of the U-NPMLE, particularly near the two end points of
interval $[0.1, 4.7]$.

%
%
%
%
%

%
%
%
%
%

%
%
%

\section{Final remarks}\label{sec6}
The estimation of the joint CDF with bivariate event time data is a
challenging problem in survival analysis. Development of
sophisticated methods for this type of problems is much
needed
{for applications}. In this paper, we develop a tensor
spline-based sieve maximum likelihood estimation method for
estimating the joint CDF with bivariate current status data. This
sieve estimation approach {reduces} the dimension of
{unknown parameter space and estimates both the joint and
marginal CDFs simultaneously}. {As a result, the proposed
method enjoys two advantages in studying bivariate event time data:
(i) it provides a unique estimate for the joint CDF, and the
numerical implementation is less demanding due to dimension
reduction; (ii) the estimation procedure automatically takes into
account the possible correlation between the two event times by
satisfying the constraints, which intuitively results in more
efficient estimation for the marginal CDFs compared to the existing
methods for estimating the marginal CDFs using only the univariate
current status data}.

Under mild regularity conditions, we also show that the proposed
spline-based sieve estimator is consistent and could converge to the
true joint CDF at a rate faster than $n^{1/3}$ if the target CDF is
smooth enough. {Both theoretical and numerical results provide
evidence that the proposed sieve MLE outperforms the conventional
NPMLE studied in the literature}. {The superior performance of the
proposed method mainly rests on the smoothness of the true bivariate
distribution function. In many applications of bivariate survival
analysis, this assumption of smoothness is reasonable and shall
motivate the use of the proposed method}.

Though the development of the proposed method is illustrated with
bivariate current status data as it algebraically simplifies the
theoretical justification, the proposed method can be
readily extended to bivariate interval censored data [\citet{Son01} and
\citet{Maa05}] as well as bivariate right censored data [\citet{Dab88} and \citet{Koo98}] {with parallel theoretical and
numerical justifications}. {It is potentially applicable in
any nonparametric estimation problem of multivariate distribution
function.}

While the consistency and rate of convergence are fully studied for
the proposed estimator, the study of its asymptotic distribution is
not accomplished. With the knowledge of asymptotic distribution of
the conventional NPMLE for current status data studied in \citet{GroWel92} and \citet{Son01}, it is for sure that the
asymptotic distribution of the proposed estimator will not be
Gaussian. Discovering the limiting distribution for the proposed
estimator remains an interesting yet a very challenging problem for
future investigation.

\section{Proofs of the theorems}\label{sec7}
For the rest of this paper, we denote $K$ as a universal positive
constant that may be different from place to place and $\mathbb
{P}_nf=\frac{1}{n}\sum_{i=1}^nf(X_i)$, the empirical process indexed by $f(X)$.

\begin{pf*}{Proof of Theorem \ref{th3.1}}
We show $\hat{\tau}_n$ is a consistent estimator by verifying the
three conditions of Theorem 5.7 in \citet{van98}.

For $(s,t)\in[l_1,u_1]\times[l_2,u_2]$, we define $\Omega$ by
\begin{eqnarray*}
\Omega&=&\bigl\{\tau(s,t)=\bigl(F(s,t),F_1(s),F_2(t)\bigr):
\\
&&\hphantom{\{} \tau\mbox{ satisfies the following conditions (a) and
(b)}\bigr\}\mbox{:}
\end{eqnarray*}

\begin{longlist}[(a)]
\item[(a)]
$F(s,t)$ is nondecreasing in both $s$ and $t$, $F_{1}(s)-F(s,t)$ is
nondecreasing in $s$ but nonincreasing in $t$, $F_{2}(t)-F(s,t)$ is
nondecreasing in $t$ but nonincreasing in $s$, and
$1-F_1(s)-F_2(t)+F(s,t)$ is nonincreasing in both~$s$ and $t$,
\item[(b)]
$F(s,t)\geq b_1$, $F_{1}(s)-F(s,t)\geq b_2$, $F_{2}(t)-F(s,t)\geq b_3$
and $1-F_{1}(s)-F_{2}(t)+F(s,t)\geq b_4$, for $b_1>0$, $b_2>0$, $b_3>0$
and $b_4>0$.
\end{longlist}

Since (C2) and (C6) hold, Lemma 0.1 in the supplemental article [\citet{WuZha}]
implies that there exist $b_1>0$, $b_2>0$,
$b_3>0$ and $b_4>0$ small enough to guarantee that $\tau_0\in\Omega$
and $\Omega_n^{\prime}\in\Omega$. We suppose $b_1$, $b_2$, $b_3$ and
$b_4$, in condition (b) above, are chosen small enough such that
$\Omega$ contains both~$\tau_0$ and $\Omega_n^{\prime}$.

Denote $\mathcal{L}=\{l(\tau)\dvtx \tau\in\Omega\}$ the class of
functions induced by the log likelihood with a single observation
$x=(s,t,\delta_1,\delta_2)$, where
\begin{eqnarray*}
l(\tau)&=&\delta_1\delta_2\log F(s,t)+
\delta_1(1-\delta_2)\log\bigl[F_1(s)-F(s,t)
\bigr]
\\
&&{}+(1-\delta_1)\delta_2\log\bigl[F_2(t)-F(s,t)
\bigr]
\\
&&{}+(1-\delta_1) (1-\delta_2)\log\bigl[1-F_1(s)-F_2(t)+F(s,t)
\bigr],
\end{eqnarray*}
with $\delta_1=1_{[T_1\leq s]},\delta_2=1_{[T_2\leq t]}$.
We denote $\mathbb{M}(\tau)=Pl(\tau)$ and $\mathbb{M}_n(\tau)=\mathbb
{P}_n(l(\tau))$.

First, we verify
$
\sup_{\tau\in\Omega}|\mathbb{M}_n(\tau)-\mathbb{M}(\tau)|\rightarrow_p 0$.

It suffices to show that $\mathcal{L}$ is a $P$-Glivenko--Cantelli, since
\[
\sup_{\tau\in\Omega}\bigl|\mathbb{M}_n(\tau)-\mathbb{M}(\tau)\bigr| =
\sup_{l(\tau)\in\mathcal{L}}\bigl|(\mathbb{P}_n-P)l(\tau)\bigr|\rightarrow_p
0.
\]

Let $A_1= \{\frac{\log F(s,t)}{\log b_1}\dvtx \tau=(F,F_1,F_2)\in
\Omega \}$, and $\mathcal{G}_1=\{1_{[l_1,s]\times[l_2,t]},
l_1\leq s\leq u_1,l_2\leq t\leq u_2\}.$ By conditions (a) and (b),
we know $0\leq\frac{\log F(s,t)}{\log b_1}\leq1$ and $\frac{\log
F(s,t)}{\log b_1}$ is nonincreasing in both $s$ and $t$. Therefore
$A_1\subseteq\overline{\operatorname{sconv}}(\mathcal{G}_1)$, the closure of the
symmetric convex hull of $\mathcal{G}_1$
[van der Vaart and Wellner (\citeyear{vanWel96})]. Hence Theorem 2.6.7 in van der Vaart and Wellner (\citeyear{vanWel96})
implies that
%
\begin{equation}
\label{g1number} N\bigl(\varepsilon, \mathcal{G}_1,
L_2(Q_{C_1,C_2})\bigr)\leq K \biggl(\frac
{1}{\varepsilon}
\biggr)^4
\end{equation}
for any probability measure $Q_{C_1,C_2}$ of $(C_1, C_2)$. By the
facts that $V(\mathcal{G}_1)=3$ and the envelop function of
$\mathcal{G}_1$ is 1. (\ref{g1number}) is followed by
\[
\log N\bigl(\varepsilon, \overline{\operatorname{sconv}}(\mathcal{G}_1),
L_2(Q_{C_1,C_2})\bigr)\leq K \biggl(\frac{1}{\varepsilon}
\biggr)^{4/3},
\]
using the result of Theorem 2.6.9 in van der Vaart and Wellner (\citeyear{vanWel96}).
Hence
%
\begin{equation}
\label{a1number} \log N\bigl(\varepsilon, A_1, L_2(Q_{C_1,C_2})\bigr)
\leq K \biggl(\frac{1}{\varepsilon
} \biggr)^{4/3}.
\end{equation}
Let
\[
A_1^{\prime}=\bigl\{\delta_1\delta_2
\log F(s,t)\dvtx \tau=(F,F_1,F_2)\in\Omega\bigr\}.
\]

Suppose the centers of $\varepsilon$-balls of $ A_1$ are $f_i,
i=1,2,\ldots,[K(\frac{1}{\varepsilon})^{4/3}]$, and then for any joint
probability measure $Q$ of $(T_1,T_2,C_1,C_2)$,
\begin{eqnarray*}
&&\|\delta_1\delta_2\log F-\delta_1
\delta_2\log b_1f_i\|^2_{L_2(Q)}\\[-2pt]
&&\qquad=Q
\biggl[\delta_1\delta_2\log b_1 \biggl(
\frac{\log F}{\log
b_1}-f_i \biggr) \biggr]^2
\\[-2pt]
&&\qquad=E \biggl[1_{[T_1<C_1,T_2<C_2]}\log b_1 \biggl(
\frac{\log
F(C_1,C_2)}{\log b_1}-f_i(C_1,C_2) \biggr)
\biggr]^2
\\[-2pt]
&&\qquad=E \biggl\{E \biggl\{ \biggl[1_{[T_1<C_1,T_2<C_2]}\log b_1
\biggl(\frac{\log F(C_1,C_2)} {
\log b_1}-f_i(C_1,C_2)
\biggr) \biggr]^2\Big|C_1,C_2 \biggr\} \biggr\}
\\[-2pt]
&&\qquad=E_{C_1,C_2} \biggl[F_0(C_1,C_2)
\log b_1 \biggl(\frac{\log
F(C_1,C_2)} {
\log b_1}-f_i(C_1,C_2)
\biggr) \biggr]^2
\\[-2pt]
&&\qquad\leq E_{C_1,C_2} \biggl[\log b_1 \biggl(
\frac{\log
F(C_1,C_2)} {\log b_1}-f_i(C_1,C_2) \biggr)
\biggr]^2
\\[-2pt]
&&\qquad=(\log b_1)^2 \biggl\llVert
\frac{\log F}{\log b_1}-f_i\biggr\rrVert^2_{L_2(Q_{C_1,C_2})}.
\end{eqnarray*}

Let $\hat{b}_1=-\log b_1$ then
$\delta_1\delta_2\log b_1f_i, i=1,2,\ldots,[K(\frac{1}{\varepsilon
})^{4/3}]$, are the centers of
$ \varepsilon\hat{b}_1$-balls of $A_1^{\prime}$. Hence by (\ref{a1number}) we
have $ \log N(\varepsilon\hat{b}_1, A_1^{\prime}, L_2(Q))\leq K (\frac
{1}{\varepsilon} )^{4/3}$,
and it follows that
\[
\int_0^1\sup_{Q}\sqrt{
\log N\bigl(\varepsilon\hat{b}_1, A_1^{\prime},L_2(Q)
\bigr)}\,d\varepsilon \leq\int_0^1\sqrt{K}
\biggl(\frac{1}{\varepsilon} \biggr)^{2/3}\,d\varepsilon <\infty.
\]

It is obvious that the envelop function of $A_1^{\prime}$ is $\hat{b}_1$,
therefore $A_1^{\prime}$ is a $P$-Donsker by Theorem 2.5.2 in van der Vaart and Wellner (\citeyear{vanWel96}).

Let
%
\begin{eqnarray*}
A_2^{\prime}&=&\bigl\{\delta_1(1-
\delta_2)\log\bigl(F_1(s)-F(s,t)\bigr)\dvtx
\tau=(F,F_1,F_2)\in \Omega\bigr\},
\\[-2pt]
A_3^{\prime}&=&\bigl\{(1-\delta_1)
\delta_2\log\bigl(F_2(t)-F(s,t)\bigr)\dvtx
\tau=(F,F_1,F_2)\in \Omega\bigr\}
\end{eqnarray*}
and
\begin{eqnarray*}
&&A_4^{\prime}=\bigl\{(1-\delta_1) (1-
\delta_2)\log\bigl(1-F_1(s)-F_2(t)-F(s,t)
\bigr)\dvtx\\[-2pt]
&&\hspace*{143pt}\qquad \tau =(F,F_1,F_2)\in\Omega\bigr\}.\vadjust{\goodbreak}
\end{eqnarray*}
Following the same arguments for showing $A_1^{\prime}$ being a $P$-Donsker,
it can be shown that $A_2^{\prime}, A_3^{\prime}$ and $A_4^{\prime}$ are all
$P$-Donsker classes. So $\mathcal{L}$ is $P$-Donsker as well.
Since $P$-Donsker is also $P$-Glivenko--Cantelli, it then follows that
$\sup_{l(\tau)\in\mathcal{L}}|(\mathbb{P}_n-P)l(\tau)|\rightarrow_p 0.$

Second, we verify
$
\mathbb{M}(\tau_0)-\mathbb{M}(\tau)\geq Kd^2(\tau_0,\tau),
$
for any $\tau\in\Omega$.

Note that
%
\begin{eqnarray*}
&&\mathbb{M}(\tau_0)-\mathbb{M}(\tau)\\
&&\qquad=P\bigl
\{l(\tau_0)-l(\tau)\bigr\}
\\
&&\qquad=P \biggl\{\delta_1\delta_2\log\frac{F_0}{F}+
\delta_1(1-\delta_2)\log\frac{F_{0,1}-F_0}{F_1-F}\\
&&\qquad\qquad{} +(1-
\delta_1)\delta_2\log\frac{F_{0,2}-F_0}{F_2-F}
\\
&&\qquad\qquad{}  +(1-\delta_1) (1-\delta_2)\log
\frac{1-F_{0,1}-F_{0,2}+F_0} {
1-F_1-F_2+F} \biggr\}
\\
&&\qquad=P_{C_1,C_2} \biggl\{F_0\log\frac{F_0}{F}
+(F_{0,1}-F_0) \log\frac{F_{0,1}-F_0}{F_1-F}\\
&&\hspace*{46pt}\qquad{} +(F_{0,2}-F_0)
\log\frac{F_{0,2}-F_0}{F_2-F}
\\
&&\hspace*{35pt}\qquad\quad{}+(1-F_{0,1}-F_{0,2}+F_0) \log
\frac{1-F_{0,1}-F_{0,2}+F_0} {
1-F_1-F_2+F} \biggr\}, %
\end{eqnarray*}
and it follows that
%
\begin{eqnarray}
\label{ch4mq2} %
\mathbb{M}(\tau_0)-
\mathbb{M}(\tau) &=&P_{C_1,C_2} \biggl\{Fm \biggl(\frac{F_0}{F} \biggr)
+(F_1-F) m \biggl(\frac{F_{0,1}-F_0}{F_1-F} \biggr)
\nonumber\\
&&\hspace*{35pt}{}+(F_2-F) m \biggl(\frac{F_{0,2}-F_0}{F_2-F} \biggr)
\\
&&\hspace*{35pt}{}+(1-F_1-F_2+F) m \biggl(\frac{1-F_{0,1}-F_{0,2}+F_0} {
1-F_1-F_2+F}
\biggr) \biggr\}, \nonumber %
\end{eqnarray}
where $m(x)=x\log(x)-x+1\geq(x-1)^2/4$ for $0\leq x\leq5.$

Since $F$ has positive upper bound,
%
\begin{eqnarray}
\label{ch4f} %
P_{C_1,C_2} \biggl\{Fm \biggl(
\frac{F_0}{F} \biggr) \biggr\} &\geq& P_{C_1,C_2} \biggl\{F \biggl(
\frac{F_0}{F}-1 \biggr)^2\Big/4 \biggr\}\ge KP_{C_1,C_2}(F_0-F)^2
\nonumber
\\[-8pt]
\\[-8pt]
\nonumber
&=&K\|F_0-F\|^2_{L_2(P_{C_1,C_2})}. %
\end{eqnarray}
%
Similarly, we can easily show that
%
\begin{eqnarray}
\label{ch4f1f} %
&& P_{C_1,C_2} \biggl
\{(F_1-F) m \biggl(\frac{F_{0,1}-F_0}{F_1-F} \biggr) \biggr\} 
\nonumber
\\[-8pt]
\\[-8pt]
\nonumber
&&\qquad\ge K\bigl\|(F_{0,1}-F_1)-(F_0-F)
\bigr\|^2_{L_2(P_{C_1,C_2})}, %
\\
\label{ch4f2f} %
&& P_{C_1,C_2} \biggl
\{(F_2-F)m \biggl(\frac{F_{0,2}-F_0}{F_2-F} \biggr) \biggr\} 
\nonumber
\\[-8pt]
\\[-8pt]
\nonumber
&&\qquad\ge
K\bigl\|(F_{0,2}-F_2)-(F_0-F)\bigr\|^2_{L_2(P_{C_1,C_2})}
\end{eqnarray}
and
%
\begin{eqnarray}
\label{ch4f1f2f} %
&&P_{C_1,C_2} \biggl
\{(1-F_1-F_2+F) m \biggl(\frac{1-F_{0,1}-F_{0,2}+F_0} {
1-F_1-F_2+F} \biggr)
\biggr\}
\nonumber
\\[-8pt]
\\[-8pt]
\nonumber
&&\qquad\geq K\bigl\|(1-F_{0,1}-F_{0,2}+F_0)-
(1-F_1-F_2+F)\bigr\|^2_{L_2(P_{C_1,C_2})}.
\end{eqnarray}

So combining (\ref{ch4f}), (\ref{ch4f1f}), (\ref{ch4f2f}) and (\ref
{ch4f1f2f}) results in
\begin{eqnarray*}
\mathbb{M}(\tau_0)-\mathbb{M}(\tau) &\geq&K \bigl(
\|F_0-F\|^2_{L_2(P_{C_1,C_2})}
\\
&&\hspace*{14pt}{}+\bigl\|(F_{0,1}-F_1)-(F_0-F)
\bigr\|^2_{L_2(P_{C_1,C_2})}
\\
&&\hspace*{14pt}{}  +\bigl\|(F_{0,2}-F_2)-(F_0-F)
\bigr\|^2_{L_2(P_{C_1,C_2})} \bigr).
\end{eqnarray*}
Let $f_1=\|F_0-F\|^2_{L_2(P_{C_1,C_2})}$,
$f_2=\|F_{0,1}-F_1\|^2_{L_2(P_{C_1})}$ and
$f_3=\|F_{0,2}-  F_2\|^2_{L_2(P_{C_2})}$. If $f_1$ is the largest
among $f_1, f_2, f_3,$ then
%
\begin{equation}
\label{ch4f1} \mathbb{M}(\tau_0)-\mathbb{M}(\tau)\geq
Kf_1\geq(K/3) (f_1+f_2+f_3).
\end{equation}
If $f_2$ is the largest, then
%
\begin{equation}
\label{ch4f2}\qquad \mathbb{M}(\tau_0)-\mathbb{M}(\tau)\geq K
\bigl[f_1+(f_2-f_1)\bigr]\geq
Kf_2\geq (K/3) (f_1+f_2+f_3).
\end{equation}
If $f_3$ is the largest, then
%
\begin{equation}
\label{ch4f3}\qquad\quad \mathbb{M}(\tau_0)-\mathbb{M}(\tau)\geq K
\bigl[f_1+(f_3-f_1)\bigr]\geq
Kf_3\geq (K/3) (f_1+f_2+f_3).
\end{equation}
Therefore, by (\ref{ch4f1}), (\ref{ch4f2}) and (\ref{ch4f3}), it
follows that
\[
\mathbb{M}(\tau_0)-\mathbb{M}(\tau)\geq Kd^2(
\tau_0,\tau).
\]

Finally, we verify
$\mathbb{M}_n(\hat{\tau}_n)-\mathbb{M}_n(\tau_0)\geq-o_p(1).$

Since (C2), (C3) and (C6) hold, Lemma 0.3 in the supplemental article [\citet{WuZha}] implies that there exists $\tau_n=(F_n, F_{n,1},
F_{n,2})$ in $\Omega_n^{\prime}$ such that for
$\tau_0=(F_0,F_{0,1},F_{0,2})$, $\|F_n-F_0\|_\infty\leq K(n^{-pv}),$
$\|F_{n,1}-F_{0,1}\|_\infty\leq K(n^{-pv})$ and
$\|F_{n,2}-F_{0,2}\|_\infty\leq K(n^{-pv}).$ Since $\hat{\tau}_n$
maximizes $\mathbb{M}_n(\tau)$ in $\Omega_n^{\prime}$,
$\mathbb{M}_n(\hat{\tau}_n)-\mathbb{M}_n(\tau_n)>0.$ Hence,
%
\begin{eqnarray}
\label{ch4mnmn} %
\mathbb{M}_n(\hat{
\tau}_n)-\mathbb{M}_n(\tau_0) &=&
\mathbb{M}_n(\hat{\tau}_n)-\mathbb{M}_n(
\tau_n)+\mathbb{M}_n(\tau_n)-
\mathbb{M}_n(\tau_0)
\nonumber\\
&\geq& \mathbb{M}_n(\tau_n)-\mathbb{M}_n(
\tau_0)=\mathbb{P}_n\bigl(l(\tau_n)\bigr)-
\mathbb{P}_n\bigl(l(\tau_0)\bigr)
\\
&=&(\mathbb{P}_n-P)\bigl\{l(\tau_n)-l(
\tau_0)\bigr\}+P\bigl\{l(\tau_n)-l(\tau_0)
\bigr\}.\nonumber %
\end{eqnarray}

Define
\begin{eqnarray*}
\mathcal{L}_{n}&=&\bigl\{l(\tau_n)\dvtx
\tau_n=(F_n,F_{n,1},F_{n,2})\in
\Omega_n^{\prime}, \|F_n-F_0
\|_{\infty}\leq K\bigl(n^{-pv}\bigr),
\\
&&\hspace*{14pt}{}\|F_{n,1}-F_{0,1}\|_{\infty}\leq K
\bigl(n^{-pv}\bigr), \|F_{n,2}-F_{0,2}
\|_{\infty}\leq K\bigl(n^{-pv}\bigr)\bigr\}.
\end{eqnarray*}

Since $(a+b+c+d)^2\leq4(a^2+b^2+c^2+d^2)$, then for any $l(\tau_n)\in
\mathcal{L}_{n}$, we have
%
\begin{eqnarray}
\label{ch4pltau} %
&&P\bigl\{l(\tau_n)-l(
\tau_0)\bigr\}^2\nonumber\\
 &&\qquad\leq 4P \biggl(\delta_1
\delta_2\log\frac{F_n}{F_0} \biggr)^2
\nonumber\\
&&\qquad\quad{}+4P \biggl(\delta_1(1-\delta_2)\log
\frac
{F_{n,1}-F_n}{F_{0,1}-F_0} \biggr)^2 +4P \biggl((1-\delta_1)
\delta_2\log\frac{F_{n,2}-F_n}{F_{0,2}-F_0} \biggr)^2
\nonumber
\\[-8pt]
\\[-8pt]
\nonumber
&&\qquad\quad{}+4P \biggl((1-\delta_1) (1-\delta_2) \log
\frac{1-F_{n,1}-F_{n,2}+F_n}{1-F_{0,1}-F_{0,2}+F_0} \biggr)^2
\\
&&\qquad \leq4P_{C_1,C_2} \biggl(\log\frac{F_n}{F_0}
\biggr)^2 +4P_{C_1,C_2} \biggl(\log\frac{F_{n,1}-F_n}{F_{0,1}-F_0}
\biggr)^2
\nonumber\\
&&\qquad\quad{}+4P_{C_1,C_2} \biggl(\log\frac
{F_{n,2}-F_n}{F_{0,2}-F_0}
\biggr)^2 +4P_{C_1,C_2} \biggl(\log\frac{1-F_{n,1}-F_{n,2}+F_n}{1-F_{0,1}-F_{0,2}+F_0}
\biggr)^2.\nonumber %
\end{eqnarray}

The facts that\vspace*{1pt} $\|F_n-F_0\|_{\infty}\leq K(n^{-pv})$ and that $F_0$ has
a positive lower bound result in
$1/2<\frac{F_n}{F_0}<2$ for large $n$. It can be easily shown that if
$1/2\leq x\leq2$, $|\log(x)|\leq K|x-1|$. Hence
$\llvert \log\frac{F_n}{F_0}\rrvert \leq K\llvert \frac{F_n}{F_0}-1\rrvert $,
and it follows that
%
\begin{equation}
\label{ch4fnf0} %
\qquad P_{C_1,C_2}\biggl\llvert \log
\frac{F_n}{F_0}\biggr\rrvert^2 \leq KP_{C_1,C_2}\biggl\llvert
\frac{F_n}{F_0}-1\biggr\rrvert^2\leq KP_{C_1,C_2}|F_n-F_0|^2
\rightarrow0. %
\end{equation}
%
Similar arguments yield
%
\begin{eqnarray}
\label{ch4fnf01} %
\qquad\quad P_{C_1,C_2}\biggl\llvert \log
\frac{F_{n,1}-F_n}{F_{0,1}-F_0}\biggr\rrvert^2 
&\leq& KP_{C_1,C_2}\bigl|(F_{n,1}-F_n)-(F_{0,1}-F_0)\bigr|^2
\rightarrow0, %
\\
\label{ch4fnf02} %
P_{C_1,C_2}\biggl\llvert \log
\frac{F_{n,2}-F_n}{F_{0,2}-F_0}\biggr\rrvert^2 &\leq& KP_{C_1,C_2}\bigl|(F_{n,2}-F_n)-(F_{0,2}-F_0)\bigr|^2
\rightarrow0 %
\end{eqnarray}
and
%
\begin{equation}
\label{ch4fnf03} %
P_{C_1,C_2}\biggl\llvert \log
\frac{1-F_{n,1}-F_{n,2}+F_n}{1-F_{0,1}
-F_{0,2}+F_0}\biggr\rrvert^2 \rightarrow0. %
\end{equation}
Combining (\ref{ch4pltau})--(\ref{ch4fnf03}) results in $P\{l(\tau_n)-l(\tau_0)\}^2\rightarrow0, \mbox{ as } n\rightarrow\infty$. Hence
%
\begin{eqnarray}
\label{ch4square3} \rho_P\bigl\{l(\tau_n)-l(
\tau_0)\bigr\}&=&\bigl\{\operatorname{var}_P\bigl[l(\tau_n)-l(
\tau_0)\bigr]\bigr\}^{1/2}
\nonumber
\\[-8pt]
\\[-8pt]
\nonumber
&\le&\bigl\{P\bigl[l(
\tau_n)-l(\tau_0)\bigr]^2\bigr
\}^{1/2}\rightarrow_{n\rightarrow\infty}0.
\end{eqnarray}

Since $\mathcal{L}$ is shown a $P$-Donsker in the first part of the
proof, Corollary~2.3.12 of van der Vaart and Wellner (\citeyear{vanWel96}) yields
that
%
\begin{equation}
\label{ch4i1} (\mathbb{P}_n-P)\bigl\{l(\tau_n)-l(
\tau_0)\bigr\}=o_p \bigl(n^{-1/2} \bigr),
\end{equation}
by the fact that both $l(\tau_n)$ and $l(\tau_0)$ are in ${\mathcal L}$
and (\ref{ch4square3}).

In addition,
\begin{eqnarray*}
\bigl\llvert P\bigl\{l(\tau_n)-l(\tau_0)\bigr\}\bigr
\rrvert \leq P\bigl|l(\tau_n)-l(\tau_0)\bigr| \leq K \bigl\{P
\bigl[l(\tau_n)-l(\tau_0)\bigr]^2 \bigr
\}^{1/2}\rightarrow_{n\rightarrow\infty}0.
\end{eqnarray*}
Therefore $P(l(\tau_n)-l(\tau_0))\geq-o(1)$ as $n\rightarrow\infty$. Hence,
\[
\mathbb{M}_n(\hat{\tau}_n)-\mathbb{M}_n(
\tau_0)\geq o_p\bigl(n^{-1/2}\bigr)-o(1)
\geq-o_p(1).
\]

This completes the proof of $d(\hat{\tau}_n, \tau_0)\rightarrow0$ in
probability.
\end{pf*}

\begin{pf*}{Proof of Theorem \ref{th3.2}}
We derive the rate of convergence by verifying the conditions of
Theorem 3.4.1 of van der Vaart and Wellner (\citeyear{vanWel96}).
To apply the theorem to this problem, we denote $M_n(\tau)={\mathbb
{M}}(\tau)=Pl(\tau)$ and $d_n(\tau_1,\tau_2)=d(\tau_1,\tau_2)$.
The maximizer of $\mathbb{M}(\tau)$ is $\tau_0=(F_0,F_{0,1},F_{0,2})$.

\begin{longlist}[(ii)]
\item[(i)] Let $\tau_n\in\Omega_n^{\prime}$ with $\tau_n$
satisfying $d(\tau_n,\tau_0)\leq K (n^{-pv})$ and
$\delta_n=n^{-pv}$. We verify that for large $n$ and any
$\delta>\delta_n$,
\[
\mathop{\operatorname{sup}}_{\delta/2 < d(\tau,\tau_n) \leq\delta, \tau\in\Omega_n'} \bigl(\mathbb{M}(\tau)-\mathbb{M}(
\tau_n) \bigr)\leq-K\delta^2.
\]

Since $d(\tau,\tau_0)\geq d(\tau,\tau_n)-d(\tau_0,\tau_n)\geq
\delta/2-K(n^{-pv})$, then for large $n$, $d(\tau,\tau_0)\geq
K\delta$. In the proof of consistency, we have already established
that $\mathbb{M}(\tau)-\mathbb{M}(\tau_0)\leq-Kd^2(\tau,\tau_0)\leq
-K\delta^2$. And as shown in the proof of consistency,
$\mathbb{M}(\tau_0)-\mathbb{M}(\tau_n)\leq Kd^2(\tau_0,\tau_n)\leq
K(n^{-2pv})$. Therefore, for large $n$,
$\mathbb{M}(\tau)-\mathbb{M}(\tau_n)=\mathbb{M}(\tau)-\mathbb{M}(\tau_0)+\mathbb{M}(\tau_0)-\mathbb{M}(\tau_n)\leq
-K\delta^2+K(n^{-2pv})=-K\delta^2$.

\item[(ii)] We {shall} find a function
$\psi(\cdot)$ such that
\[
E \Bigl\{\mathop{\operatorname{sup}}_{{\delta/2 < d(\tau,\tau_n) \leq\delta, \tau\in\Omega_n'}
} \mathbb{G}_n(
\tau-\tau_n) \Bigr\}\leq K\frac{\psi(\delta)}{\sqrt{n}}
\]
and $\delta\rightarrow\psi(\delta)/\delta^{\alpha}$ is decreasing on
$\delta$, for some $\alpha<2$, and for $r_n\leq\delta_n^{-1}$, it satisfies
\[
r^2_n\psi (1/r_n )\le K\sqrt{n} \qquad
\mbox{for every } n.
\]

Let
\[
\mathcal{L}_{n,\delta}=\bigl\{l(\tau)-l(\tau_n)\dvtx \tau\in
\Omega_n^{\prime} \mbox{ and } \delta/2<d(\tau,
\tau_n)\leq\delta\bigr\}.
\]
First, we evaluate the bracketing number of $\mathcal{L}_{n,\delta}$.
\end{longlist}

Let $\mathcal{F}_n=\{F\dvtx \tau=(F,F_{1},F_{2})\in\Omega_n^{\prime}, \delta
/2\leq d(\tau,\tau_n)\leq\delta\}$,
$\mathcal{F}_{n,1}=\{F_{1}\dvtx \tau=(F,F_{1},F_{2})\in\Omega_n^{\prime}, \delta
/2\leq d(\tau,\tau_n)\leq\delta\}$ and
$\mathcal{F}_{n,2}=\{F_{2}\dvtx \tau=(F,F_{1},F_{2})\in\Omega_n^{\prime},\break  \delta
/2\leq d(\tau,\tau_n)\leq\delta\}$.\vadjust{\goodbreak}

Denote $\tau_n=(F_n,F_{n,1},F_{n,2})$. Lemma 0.5 in the supplemental
article [\citet{WuZha}] implies that there exist
$\varepsilon$-brackets $[D^{L}_i, D^{U}_i],
i=1,2,\ldots,\break [(\delta/\varepsilon)^{Kp_nq_n}]$ to cover
$\mathcal{F}_n-F_n$. Moreover, Lemma 0.6 in the supplemental article
[\citet{WuZha}] implies there exist $\varepsilon$-brackets
$[D^{(1),L}_j, D^{(1),U}_j],
j=1,2,\ldots,[(\delta/\varepsilon)^{Kp_n}]$, to cover
$\mathcal{F}_{n,1}-F_{n,1}$, and there exist $\varepsilon$-brackets
$[D^{(2),L}_k,\break  D^{(2),U}_k],
k=1,2,\ldots,[(\delta/\varepsilon)^{Kq_n}]$, to cover
$\mathcal{F}_{n,2}-F_{n,2}$.

Denote $F^L_i\equiv D^L_i+F_n$, $F^U_i\equiv D^U_i+F_n$,
$F^{(1),L}_j\equiv D^{(1),L}_j+F_{n,1}$, $F^{(1),U}_j\equiv
D^{(1),U}_j+F_{n,1}$, $F^{(2),L}_k\equiv D^{(2),L}_k+F_{n,2}$ and
$F^{(2),U}_k\equiv D^{(2),U}_k+F_{n,2}$. Let
\begin{eqnarray*}
l^U_{i,j,k}&=&\delta_1\delta_2\log
F^{U}_i +\delta_1(1-\delta_2)
\log\bigl(F^{(1),U}_j-F^{L}_i\bigr)
\\
&&{}+(1-\delta_1)\delta_2\log\bigl(F^{(2),U}_k-F^{L}_i
\bigr)
\\
&&{}+(1-\delta_1) (1-\delta_2)\log\bigl(1-F^{(1),L}_j-F^{(2),L}_k+F^{U}_i
\bigr)
\end{eqnarray*}
and
\begin{eqnarray*}
l^L_{i,j,k}&=&\delta_1\delta_2\log
F^{L}_i +\delta_1(1-\delta_2)
\log\bigl(F^{(1),L}_j-F^{U}_i\bigr)
\\
&&{}+(1-\delta_1)\delta_2\log\bigl(F^{(2),L}_k-F^{U}_i
\bigr)
\\
&&{}+(1-\delta_1) (1-\delta_2)\log\bigl(1-F^{(1),U}_j-F^{(2),U}_k+F^{L}_i
\bigr).
\end{eqnarray*}
Then for any $l(\tau)\in\{\mathcal{L}_{n,\delta}+l(\tau_n)$\}, there
exist $i,j,k,$
for $i=1,2,\ldots,\break[(\delta/\varepsilon)^{Kp_nq_n}]$,
$j=1,2,\ldots,[(\delta/\varepsilon)^{Kp_n}]$ and
$k=1,2,\ldots,[(\delta/\varepsilon)^{Kq_n}]$, such that
$l^L_{i,j,k}\leq l(\tau)\leq l^U_{i,j,k}$ and the number of brackets
$[l^L_{i,j,k}, l^U_{i,j,k}]^{\prime}s$ is bounded by
$(\delta/\varepsilon)^{Kp_nq_n}\cdot(\delta/\varepsilon)^{Kp_n}\cdot(\delta
/\varepsilon)^{Kq_n}.$

Note that
\begin{eqnarray*}
&&\bigl\Vert  l^U_{i,j,k}-l^L_{i,j,k}
\bigr\Vert _\infty\\
&&\qquad\leq 
\biggl\llVert \log
\frac{F^{U}_i}{F^{L}_i}\biggr\rrVert_\infty +\biggl\llVert \log
\frac{F^{(1),U}_j-F^{L}_i} {
F^{(1),L}_j-F^{U}_i}\biggr\rrVert_\infty
\\
&&\qquad\quad{}+\biggl\llVert \log\frac{F^{(2),U}_k-F^{L}_i} {
F^{(2),L}_k-F^{U}_i}\biggr\rrVert_\infty +\biggl
\llVert \log\frac{1-F^{(1),L}_j-F^{(2),L}_j
+F^{U}_i} {
1-F^{(1),U}_j-F^{(2),U}_j
+F^{L}_i}\biggr\rrVert_\infty.
\end{eqnarray*}

Since for any $\tau\in\Omega_n'$, $F$ has a positive lower bound,
then for a small $\varepsilon$, $F^{L}_i$ can be made to have a
positive lower bound as well.\vspace*{-1pt} Combining with the fact that
$F^{U}_i(s,t)$ is
close to $F^{L}_i(s,t)$ guarantees that $0\leq\frac
{F^{U}_i}{F^{L}_i}-1\leq1$
for $i=1,2,\ldots,[(\delta/\varepsilon)^{Kp_nq_n}]$.
Note that by $\log x\leq(x-1)$ for $0\leq(x-1)\leq1$, therefore $\log
\frac{F^{U}_i}{F^{L}_i}\leq\frac{F^{U}_i} {
F^{L}_i}-1.$

Hence,
\[
\label{c4s31} \biggl\llVert \log\frac{F^{U}_i}{F^{L}_i}\biggr\rrVert_\infty
\leq\biggl\llVert \frac{F^{U}_i} {
F^{L}_i}-1\biggr\rrVert_\infty \leq\biggl
\llVert \frac{1}{F^{L}_i}\bigl(F^{U}_i-F^{L}_i
\bigr)\biggr\rrVert_\infty \leq K\bigl\llVert F^{U}_i-F^{L}_i
\bigr\rrVert_\infty\leq K\varepsilon.
\]


Similarly, by the definition of $\Omega_n'$, we can easily show that
\[
\biggl\llVert \log\frac{F^{(1),U}_j-F^{L}_i} {
F^{(1),L}_j-F^{U}_i}\biggr\rrVert_\infty 
\leq K\varepsilon, \qquad
\ 
\biggl\llVert \log\frac{F^{(2),U}_k-F^{L}_i} {
F^{(2),L}_k-F^{U}_i}\biggr\rrVert_\infty
\leq K\varepsilon
\]
and
\[
\biggl\llVert \log\frac{1-F^{(1),L}_j-F^{(2),L}_j
+F^{U}_i} {
1-F^{(1),U}_j-F^{(2),U}_j
+F^{L}_i}\biggr\rrVert_\infty \leq K
\varepsilon.
\]

Hence, the fact that $L_2$-norm is bounded by $L_\infty$-norm results in
%
\begin{equation}
\label{c4s3bn} N_{[\,]}\bigl\{\varepsilon, \mathcal{L}_{n,\delta},
L_2(P)\bigr\}\leq N_{[\, ]}\bigl\{\varepsilon,
\mathcal{L}_{n,\delta}, \Vert \cdot \Vert _\infty\bigr\}
\leq(\delta/\varepsilon)^{Kp_nq_n}.
\end{equation}

{Next}, we show that $P\{l(\tau)-l(\tau_n)\}^2\leq K\delta^2$ for
any $l(\tau)-l(\tau_n)\in\mathcal{L}_{n,\delta}$. Since for any
$\tau=(F,F_1,F_2)$ with $d(\tau,\tau_n)<\delta$,
$\|F-F_n\|_{L_2(P_{C_1,C_2})}\leq d(\tau,\tau_n)\leq\delta$. Then
with (C1), (C3) and (C5), Lemma 0.7 in the supplemental article [\citet{WuZha}]
implies that for a small $\delta>0$ and a
sufficiently large $n$, $F$ and $F_n$ are both very close to $F_0$
at every point in $[l_1,u_1]\times[l_2,u_2]$. Therefore, $F$ and
$F_n$ are very close to each other at every point in
$[l_1,u_1]\times[l_2,u_2]$. Then the fact that $F_n$ has a positive
lower bound results in $1/2<\frac{F}{F_n}<2$.
Hence $\llvert \log\frac{F}{F_n}\rrvert \leq K\llvert \frac{F}{F_n}-1\rrvert $, and it follows that
\[
\label{ch4r1} %
P_{C_1,C_2}\biggl\llvert \log
\frac{F}{F_n}\biggr\rrvert^2 \leq KP_{C_1,C_2}\biggl\llvert
\frac{F}{F_n}-1\biggr\rrvert^2 
\leq
KP_{C_1,C_2}\llvert F-F_n\rrvert^2\leq K
\delta^2. %
\]


Again by the definition of $\Omega_n'$, we can similarly show that,
given a small $\delta>0$,
when $n$ is large enough, the following inequalities are true:
\[
P_{C_1,C_2}\biggl\llvert \log\frac{F_{1}-F}{F_{n,1}-F_n}\biggr\rrvert^2
\leq K\delta^2, \qquad
P_{C_1,C_2}\biggl\llvert \log\frac{F_{2}-F}{F_{n,2}-F_n}\biggr\rrvert^2
\leq K\delta^2
\]
and
\[
\label{ch4r4} %
P_{C_1,C_2}\biggl\llvert \log
\frac{1-F_{1}-F_{2}+F}{1-F_{n,1}
-F_{n,2}+F_n}\biggr\rrvert^2\leq K\delta^2.
\]

Hence for any $l(\tau)-l(\tau_n)\in\mathcal{L}_{n,\delta}$,
{it is true that}
$P\{l(\tau)-l(\tau_n)\}^2\leq K\delta^2$.
It is obvious that $\mathcal{L}_{n,\delta}$ is uniformly bounded by the
structure of the log likelihood.
Lemma 3.4.2 of van der Vaart and Wellner (\citeyear{vanWel96}) indicates that
\[
E_P\Vert \mathbb{G}_n\Vert _{\mathcal{L}_{n,\delta}} \leq
K \tilde{J}_{[\,]}\bigl\{\delta, \mathcal{L}_{n,\delta},
L_2(P)\bigr\} \biggl[1+\frac{\tilde{J}_{[\,]}\{\delta, \mathcal
{L}_{n,\delta}, L_2(P)\}}{\delta^2\sqrt{n}} \biggr],
\]
where
\begin{eqnarray*}
\tilde{J}_{[\,]}\bigl\{\delta, \mathcal{L}_{n,\delta},
L_2(P)\bigr\}= \int_0^{\delta}\sqrt{1+
\log N_{[\,]}\bigl\{\varepsilon, \mathcal {L}_{n,\delta},
L_2(P)\bigr\}}\,d\varepsilon\leq K(p_nq_n)^{1/2}
\delta, 
\end{eqnarray*}
by (\ref{c4s3bn}). This gives $ \psi(\delta)=(p_nq_n)^{1/2}\delta
+(p_nq_n)/(n^{1/2})$. It is easy to see that $\psi(\delta)/\delta$ is a
decreasing function of $\delta$. Note that for $p_n=q_n=n^v$,
\begin{eqnarray*}
n^{2pv}\psi\bigl(1/n^{pv}\bigr)&=n^{2pv}n^vn^{-pv}+n^{2pv}n^{2v}n^{-1/2}=n^{1/2}
\bigl\{ n^{pv+v-1/2}+n^{2pv+2v-1}\bigr\}.
\end{eqnarray*}
Therefore, if $pv\leq(1-2v)/2$, $n^{2pv}\psi(1/n^{pv})\leq
2n^{1/2}$. Moreover,
$n^{1-2v}\times \psi(1/\break n^{(1-2v)/2})=2n^{1/2}$. This
implies if $r_n=n^{\min\{pv,(1-2v)/2\}}$, then $r_n \leq
\delta_n^{-1}$ and $r_n^2\psi(1/r_n)\leq Kn^{1/2}$.

It is obvious that $\mathbb{M}(\hat{\tau}_n)-\mathbb{M}(\tau_n)\ge
0$ and $d(\hat{\tau}_n,\tau_n)\leq
d(\hat{\tau}_n,\tau_0)+d(\tau_0,\tau_n)\rightarrow0$ in probability.
Therefore, it follows by Theorem 3.4.1 in van der Vaart and Wellner (\citeyear{vanWel96}) that $r_nd(\hat{\tau}_n,\tau_n)=O_p(1)$. Hence, by
$d(\tau_n,\tau_0)\leq K (n^{-pv})$
\[
r_nd(\hat{\tau}_n,\tau_0)\leq
r_nd(\hat{\tau}_n,\tau_n)+r_nd(
\tau_n,\tau_0)=O_p(1).
\]
\upqed\end{pf*}

\section*{Acknowledgments}
We owe thanks to the Associate Editor and two anonymous referees for their helpful and constructive
comments and suggestions that helped improve the manuscript from an early version.

\begin{supplement}[id=suppA]
\stitle{Technical lemmas}
\slink[doi]{10.1214/12-AOS1016SUPP} 
\sdatatype{.pdf}
\sfilename{AOS1016\_supp.pdf}
\sdescription{This supplemental material contains some technical
lemmas including their proofs that are imperative for the proofs of
Theorems \ref{th3.1} and~\ref{th3.2}.}
\end{supplement}




\printaddresses

\end{document}